\documentclass [11pt]{amsart} 
\usepackage{epic,eepic,latexsym, amssymb, amscd, amsfonts}

\input xy
\xyoption {all}
\def \uv {{\mathfrak M_{v}}}
\def \uw {{\mathfrak M_{w}}}
\def \fm {{{\bf R}{\mathcal S}}}
\def \rp {{{\bf R}}}

\newtheorem{theorem}{Theorem} 
\newtheorem {lemma}{Lemma}     
\newtheorem{conjecture}{Conjecture} 
\newtheorem {corollary}{Corollary} 

\theoremstyle{definition}
\newtheorem{remark}{Remark} 
\newtheorem{example}{Example} 
\newtheorem{convention}{Convention} 

\theoremstyle {definition}

\begin{document}

\baselineskip=17pt

\title[Bundles of generalized theta functions]{Bundles of generalized theta functions over abelian surfaces}
\author {Dragos Oprea}
\address {Department of Mathematics}
\address {University of California San Diego}
\email {doprea@math.ucsd.edu}
\date{} 

\begin {abstract} 
We study the Verlinde bundles of generalized theta functions constructed from moduli spaces of sheaves over abelian surfaces. In degree $0$, the splitting type of these bundles is expressed in terms of indecomposable semihomogeneous factors. Furthermore, Fourier-Mukai symmetries of the Verlinde bundles are found, consistently with strange duality. Along the way, a transformation formula for the theta bundles is derived, extending a theorem of Dr\'ezet-Narasimhan from curves to abelian surfaces. 
\end{abstract} 
\maketitle

\section{Introduction}

\subsection {Overview} In this paper, we put forward an analogy between aspects the strange duality proposal for curves and abelian surfaces. Such an analogy is by no means obvious:  we study moduli spaces whose geometries are very different. It is therefore surprising that the emerging pictures share common features, some which we point out below. 

In short, our main results are:
\begin {itemize}
\item [(i)]  we discuss how the theta line bundles over moduli spaces of sheaves on abelian surfaces depend on the choice of reference sheaf; the case of curves was considered in \cite {DN}; 
\item [(ii)] in degree $0$, we determine the splitting type of the Verlinde bundles (defined in Section \ref{verbdl}) in terms of indecomposable factors. In particular, we determine the action of a certain group of torsion points on the space of generalized theta functions; the curve case was solved in \cite {wirtinger};
\item [(iii)] via (i), we recast the abelian surface strange duality conjecture as a {\it specific} isomorphism between the Verlinde bundle and its Fourier-Mukai transform, as in \cite {Po}. Furthermore, using (ii), we confirm that the relevant bundles of generalized theta functions are indeed {\it abstractly} isomorphic, in degree $0$. 

\end {itemize}
As already mentioned, the results above parallel the case of curves. This is indeed one of the main points, giving credence to the strange duality conjecture in the abelian surface context. The proofs however require new ideas, and the results may be of interest independently of strange duality. 

We now detail the discussion. 

\subsection {Moduli of sheaves and their Albanese maps} To set the stage, consider a complex polarized abelian surface $(A, \Theta)$, and fix the Mukai type $v$ of sheaves $E\to A$, so that $$v=\text{ch}(E).$$ The Mukai vectors used in this paper are of the form $$v=(r, k\Theta, \chi),$$ and the following assumption will be made throughout: 

\begin {itemize}
\item [({\bf A.1})] the vector $v=(r, k\Theta, \chi)$ is primitive of positive rank, the polarization $\Theta$ is generic\footnote{this ensures that the moduli spaces we consider consist of stable sheaves only}, and furthermore the Mukai self pairing $$d_v:=\frac{1}{2}\langle v, v\rangle:=k^2\cdot \frac{\Theta^2}{2}-r\chi$$ is an odd positive integer. 
\end {itemize}

We consider the moduli space $\uv$ of $\Theta$-(semi)stable sheaves of topological type $v$. The moduli space comes equipped with the Albanese morphism $$\alpha_v=(\alpha^{+}, \alpha^{-}): \uv\to A\times \widehat A,$$ which, up to the choice of a reference sheaf $E_0$, takes sheaves $E$ to their determinant and determinant of the Fourier-Mukai transform $$\alpha_v(E)=(\det \fm(E)\otimes \det \fm (E_0)^{\vee}, \det E\otimes \det E_0^{\vee}).$$ Here $$\fm:\mathbf D(A)\to \mathbf D(\widehat A)$$ denotes the Fourier-Mukai transform. The Albanese fiber will be denoted by $K_v$, thus parametrizing semistable sheaves with fixed determinant and fixed determinant of their Fourier-Mukai transform. It is known that $K_v$ is a holomorphic symplectic manifold of dimension $2d_v-2$, deformation equivalent to the generalized Kummer variety of the same dimension \cite {Y}.

\subsection {Theta bundles} We consider the natural theta line bundles over the above moduli spaces. 
Assume that $w$ is a Mukai vector orthogonal to $v$ in the sense that in $K$-theory we have $$\chi(v\cdot w)=0.$$ Pick a complex $F\to A$ representing the Mukai vector $w$, and following \cite {jli}, \cite{LP}, construct the Fourier-Mukai transform of $F$ with kernel the universal sheaf $\mathcal E\to \uv\times A$:\footnote {or by descent from the Quot scheme in the absence of the universal sheaf} $$\Theta_w=\det \rp p_{!} (\mathcal E\otimes^{\mathbf L} q^{\star}F)^{-1}.$$ {\it Generalized theta functions} are sections of $\Theta_w$ over either one of the moduli spaces $K_v$ or $\uv$ considered above. 

\subsection {Relating different theta bundles} Over the moduli space $\uv$, the notation $\Theta_w$ is slightly imprecise, since the bundle $\Theta_F\to \uv$ may depend on the choice of representative $F$. The following result, paralleling the Dr\'ezet-Narasimhan theorem for bundles over curves \cite {DN}, controls this imprecision: 

\begin {theorem}\label{tpt} If $F_1$ and $F_2$ have the same Mukai vector orthogonal to $v$ we have $$\Theta_{F_1}=\Theta_{F_2}\otimes \left((-1)\circ \alpha^{+}\right)^{\star}\left(\det F_1\otimes \det F_2^{-1}\right)\otimes \left((-1)\circ \alpha^{-}\right)^{\star} \left(\det \fm (F_1)\otimes \det \fm (F_2)^{-1}\right).$$ 
\end {theorem} 
\noindent By contrast, the notation $\Theta_w\to K_v$ is unambiguous. 

\subsection {Verlinde numbers} \label {vnb} The holomorphic Euler characteristics of the line bundles $\Theta_w\to K_v$ are calculated in \cite {abelian}: \begin{equation}\label{eulerc}\chi(K_v, \Theta_{w})=\frac{d_v^2}{d_v+d_w}\binom {d_v+d_w}{d_v}.\end{equation} In order to use these numerics for the study of generalized theta functions, one needs to prove that $$h^0(K_v, \Theta_w)=\chi(K_v, \Theta_w).$$ This occurs for instance when $\Theta_w$ carries no higher cohomology over $K_v$, or alternatively over some smooth birational model of $K_v$. This is the case when \begin{itemize} 
\item [({\bf A.2})] \hskip.05in $\Theta_w$ belongs to the movable cone of $K_v$. \end{itemize} 
By \cite {BMOY}, the above requirement is satisfied, symmetrically in $v$ and $w$, for primitive Mukai vectors of the form\vskip.1in
\begin {itemize} 
\item [({\bf A.2})$'$] \hskip.05in $v=(r, k\Theta, \chi),\,\,\,w=(r', k'\Theta, \chi') \text{ with } k, k'>0 \text{ and }\chi, \chi'<0.$ \vskip.02in \noindent In addition, if $(k, \chi) \text { or } (k', \chi')=(1, -1)$, then $(A, \Theta)$ is not a product of elliptic curves. 
\end {itemize} 

\subsection {Bundles of generalized theta functions} \label{verbdl} When assumptions $({\bf A1})-({\bf A2})$ are satisfied, we push forward $\Theta_{w}\to \mathfrak M_v$ via the Albanese morphism $\alpha_v$. We obtain in this fashion the Verlinde bundle of generalized theta functions 
$$ \mathbf E(v,w)=(\alpha_v)_{\star} \Theta_w$$ over the abelian four-fold $A\times \widehat A$. By \eqref{eulerc} its rank equals $$\frac{d_v^2}{d_v+d_w}\binom {d_v+d_w}{d_v}.$$ 
From the discussion above, it is clear that the Verlinde bundle is well-defined only up to translation. In Sections \ref{setup} and \ref{ptb}, in particular equations \eqref{normalb} and \eqref{conv1}, we will fix the normalization of the Albanese morphism  and of the theta bundle $\Theta_w\to \mathfrak M_v$ used throughout this paper, thus pinning down the Verlinde bundle unambigously.

In the context of curves, the Verlinde bundles were introduced and studied in \cite {Po}, and were further analyzed in \cite {elliptic}, \cite {wirtinger}. For abelian surfaces, we point out how the techniques of \cite {elliptic} and \cite {wirtinger} need to be changed to the new setup we consider. 

\subsection {Semihomogeneous bundles} Assume $(A, \Theta)$ is a {\it principally} polarized surface \footnote{Most of our results also hold without change for nonprincipal polarizations. We focus on principal polarizations to keep the numerics simple.}, such that the line bundle $\Theta$ is symmetric $$(-1)^{\star}\Theta=\Theta.$$ We showed in \cite {wirtinger}, in the context of a calculation for curves, that for any pair of coprime odd integers $(a, b)$ there exists a unique symmetric semihomogeneous bundle $\mathbf W_{a, b}$ over $A$ such that $$\text { rank }\mathbf W_{a, b}=a^2 \text { and } \det \mathbf W_{a,b}=\Theta^{ab}.$$ The bundles ${\mathbf W}_{a, b}$ are the higher dimensional analogues of the Atiyah bundles over elliptic curves. Recall that semihomogeneity is the requirement that all translations of $\mathbf W_{a, b}$ by $x\in A$ are of the form $$t_x^{\star} \mathbf W_{a, b}=\mathbf W_{a, b} \otimes y$$ for some line bundle $y$ over $A$. The bundles ${\mathbf W}_{a, b}$ split after pullback \begin{equation}\label{split} a^{\star} {\mathbf W}_{a, b}=\bigoplus_{i=1}^{a^2} \Theta^{ab}.\end{equation}

We may also consider semihomogeneous bundles over the dual abelian variety $(\widehat A, -\widehat \Theta)$, where $-\widehat \Theta=\det \mathbf R\mathcal S(\Theta)^{-1}$ is the dual polarization. We will use the notation $\mathbf W_{a, b}^\dagger$ for the corresponding bundles of rank $a^2$ and determinant $\widehat{\Theta}^{ab}$. 

\subsection {Action of torsion points on generalized theta functions} \label{tp} We will now consider the case $c_1(v)=0$ {\it i.e.} $$v=(r, 0, \chi),$$ with $r$ and $\chi$ odd coprime integers, cf. {\bf (A.1)}. The dual vector $w$ must have the form $$w=(rh, k\Theta, -\chi h)$$ for some integers $h$ and $k$. 

We consider the action of $(x,y)\in A \times \widehat A$ on the moduli space $\uv$ given by $$(x,y): E\mapsto t_x^{\star} E\otimes y.$$ The action is seen to leave $K_v$ invariant provided $\chi x=0,\,\, ry=0.$ Let us write $$\gcd(\chi, k)=a, \,\,\gcd (r, k)=b.$$ If the stronger condition 
$$a x=0,\,\, by=0$$ is satisfied, the action lifts to the line bundle $\Theta_w\to K_v$, see Section \ref{s3}. We will assume this is the case. Write $$\chi(K_v, \Theta_w)=\sum_{i} (-1)^i H^i(K_v, \Theta_w)$$ for the signed sum of cohomologies of the theta bundle. We show  
\begin {theorem}\label{trcal} If $\zeta=(x,y)\in A[a]\times \widehat A[b]$ has order $\delta$, then the trace of the action of $\zeta$ is given by $$\text {Trace }(\zeta, \chi(K_v, \Theta_w))=\frac{d_v^2}{d_v+d_w}\binom{d_v/\delta+d_w/\delta}{d_v/\delta}.$$ 
\end {theorem} 

The Theorem should hold for arbitrary $c_1(v)$, but we are unable to prove this here.

\subsection {Explicit expressions for the Verlinde bundles}\label{19} We keep the same setup as in Subsection \ref{tp}. In addition, we make the assumption $({\bf A.2})$ requiring that $\Theta_w$ belong to the movable cone of $K_v$. \footnote{The strict inequalities in {\bf (A.2)$'$} are not fulfilled for the numerics we consider. However, $({\bf A.2})$ is achieved for sufficiently large slope $\mu(w)> \mu^+.$ By the argument in Section 6 of \cite {BMOY}, an explicit bound is given by $$\mu^+=\frac{\sqrt{-\chi}}{r} \cdot \frac{r-1}{\sqrt{r-2}},$$ whenever $\chi\neq -1$ and $r>2$.}

We give an explicit expression for the Verlinde bundles. The situation is easily understood when  $$\gcd (r, k)=\gcd (\chi, k)=1.$$ Then $$\mathbf E(v, w)=\bigoplus \left(\mathbf W_{-\chi, k}\boxtimes \mathbf W_{r, -k}^{\dagger} \otimes \mathcal P^{-h}\right).$$ This much can easily be derived from the representation theory of Heisenberg groups. 

The difficulty of the calculation lies however in the case when the integers $k$ and $r\chi$ are not coprime. Representation theory only gives $\mathbf E(v,w)$ up to torsion line bundles over $A\times \widehat A$ of orders  dividing $(a,b).$ We will prove the following:

\begin {theorem} \label{thm2}We have $$\mathbf E(v, w)=\bigoplus_{\zeta}\left(\mathbf W_{-\frac{\chi}{a}, \frac{k}{a}}\boxtimes \mathbf W^{\dagger}_{\frac{r}{b}, -\frac {k}{b}}\otimes \mathcal P^{-h}\right)\otimes \ell_{\zeta}^{\,\oplus \mathsf m_{\zeta}}.$$ The sum is taken over torsion line bundles $\zeta$ over $A\times \widehat A$ of orders dividing $(a, b).$ A line bundle $\zeta$ of order exactly $\omega$ comes with multiplicity $$\mathsf m_{\zeta}=\frac{1}{d_v+d_w} \sum_{\delta|ab} \frac {\delta^4}{(ab)^2} \left\{\frac{ab/\omega}{\delta}\right\} \binom{d_v/\delta+d_w/\delta}{d_v/\delta}.$$
\end {theorem} 

The line bundles $\ell_\zeta\to A\times \widehat A$ in the sum are roots of $\zeta$ of order $\left(-\frac{\chi}{a}, \frac{r}{b}\right):$ $$\left(-\frac{\chi}{a}, \frac{r}{b}\right)\ell_{\zeta}=\zeta.$$ For each $\zeta$, only one such root $\ell_{\zeta}$ is chosen. We will see that the choice of $\ell_{\zeta}$ does not affect the expressions involved. 

The Jordan totient $\left\{\,\,\,\right\}$ appearing above is defined in terms of prime factorization. Specifically, for any integer $h\geq 2$, we decompose $$h=p_1^{a_1}\ldots p_n^{a_n}$$ into powers of primes. We set $$\left\{\frac{\lambda}{h}\right\}=\begin {cases} 0 & \text {if } p_1^{a_1-1}\ldots p_n^{a_n-1} \text { does not divide } \lambda,\\ {\prod}_{i=1}^{n} \left(\epsilon_{i}-\frac{1}{p_i^{4}}\right) & \text {otherwise},\end {cases}$$ where $$\epsilon_{i}=\begin {cases} 1 & \text{ if } p_{i}^{a_{i}} |\lambda,\\ 0 & \text {otherwise.} \end {cases}$$ The Jordan totient is set to $1$ if $h=1$. 
 
\subsection {Fourier-Mukai symmetries} We furthermore consider the interaction of the Verlinde bundles with the Fourier-Mukai transform. Assuming both vectors $v$ and $w$ satisfy $({\bf A.1})-({\bf A.2})$, we show 

 \begin {theorem} \label{tttt1} When $c_1(v)$ and $c_1(w)$ are divisible by their ranks $r$ and $r'$, there is an isomorphism \begin {equation}\label{iso} \mathbf E(v, w)^{\vee}\cong \widehat{\mathbf E(w, v)}.\end{equation}
 \end {theorem}

The isomorphism \eqref{iso} is obtained by direct comparison of both sides, using Theorem \ref{thm2}. The same result should be true for any vectors $v$ and $w$ satisfying $({\bf A.1})-({\bf A.2})$. 

\subsection {Strange duality} 

The above symmetry of the Verlinde bundles is related to the strange duality conjecture\footnote {for surfaces, strange duality phenomena have first been studied by Le Potier \cite {lepotier2}}. This was observed in the case of curves by Popa \cite {Po}. 
  
Considering the fibers of \eqref{iso} over the origin, we obtain isomorphic spaces \begin{equation}\label{sdi}H^0(K_v, \Theta_w)^{\vee}\cong H^0(\uw, \Theta_v).\end{equation} As stated, this is merely saying that the dimensions of both vector spaces agree. However, there is a geometrically induced map, called {\it strange duality}, which conjecturally yields the isomorphism above. Even stronger, various strange duality maps can be packaged into an explicit bundle morphism, constructed as a corollary of Theorem \ref{tpt} \begin{equation}\label{sdm}\mathsf {SD}: \mathbf E(v, w)^{\vee}\to \widehat{\mathbf E(w, v)}.\end{equation} In many cases, ${\mathsf {SD}}$ should provide a {\it specific} geometric isomorphism of bundles, as predicted by Theorem \ref{tttt1}. This is proven for generic abelian surfaces in \cite {BMOY} for an infinite class of topological types; the general case is however still open. 

\subsection {Comparison to the case of curves} To end this introduction, let us remark that a similar picture emerges in the case of curves, cf. \cite{wirtinger}. Indeed, let $C$ be smooth of genus $g$, and consider the vectors $$v=r\left[\mathcal O_C\right], w=k\left[\kappa\right],$$ for a Theta characteristic $\kappa$.  The bundles of rank $r$, level $k$ generalized theta functions $$\mathbf E_{r, k}={\det}_{\star}\left( \Theta_{\kappa}^k\right)\to \text{Jac}(C)$$ can be obtained pushing forward the pluri-theta bundles $\Theta_w=\Theta_{\kappa}^k\to \mathfrak M_v$ via the determinant/Albanese map $$\det:\mathfrak M_v\to \text{Jac}(C).$$ The pushforwards $\mathbf E_{r, k}$ take the form $$\mathbf E_{r, k}=\bigoplus_{\zeta} \mathbf W_{\frac{r}{a}, \frac{k}{a}}\otimes \ell_{\zeta}^{\,\oplus \mathsf m_\zeta(r, k)}.$$ The torsion line bundles $\zeta$ have orders dividing $a=\gcd(r, k)$, and the $\ell_{\zeta}$s are $\frac{r}{a}$-roots of $\zeta$. The multiplicities $\mathsf m_{\zeta}(r, k)$ of $\ell_{\zeta}$ are explicit. Furthermore, the symmetry $$\mathbf E_{r, k}^{\vee} \cong \widehat {\mathbf E}_{k,r}$$ holds true, and is a manifestation of strange duality.\vskip.1in

\subsection {Outline} The paper is organized as follows. The next section discusses preliminary results about theta bundles and their behavior under \'etale pullbacks; in particular Theorem \ref{ct1} is proved there. The third section is the heart of the paper and contains the computation confirming Theorem \ref{trcal}. This is the most involved of our calculations, and Theorems \ref{thm2} and \ref{tttt1} follow from it. 

\subsection {Acknowledgements} The author was supported by the NSF through grant DMS 1150675 and by a Hellman Fellowship. Some of our results were announced in Fall 2011 under more restrictive technical assumptions, which meantime have been relaxed. 

\section {Theta bundles over the moduli of sheaves}\label{sec2}

This section collects various observations about theta bundles. The main result here is Theorem \ref{ct1}, which describes how the theta bundles depend on choices. Over curves, a similar statement was made in \cite {DN}, and proved by entirely different methods. As a corollary of the theorem, we construct the strange duality map between the Verlinde bundles. 

\subsection {Setup} \label{setup}Let $(A, \Theta)$ be a principally polarized abelian surface, with $\Theta$ a symmetric line bundle $$(-1)^{\star}\Theta=\Theta.$$ Throughout the paper, we will use the Fourier-Mukai transform $$\fm: \mathbf D(A)\to \mathbf D(\widehat A),\,\,\,\, \fm (E)=\mathbf Rp_{!}(\mathcal P\otimes q^{\star} E),$$ where $\mathcal P$ is the normalized Poincar\'e bundle over $A\times \widehat A$. If $E$ satisfies the index theorem \cite {mukai}, we often write $\widehat E$ for the sheaf representing $\fm(E)$, up to shift.  We furthermore recall the following identities \cite {mukai}: $$\fm (t_{x}^{\star} E)=\fm (E)\otimes\mathcal P_{-x},$$ $$\fm (E\otimes y)=t_{y}^{\star}\, \fm (E),$$ for $x\in A, \,y\in \widehat A$. We set $$\widehat \Theta=\det \fm (\Theta),$$ so that $-\widehat \Theta$ is the polarization on the dual abelian variety $\widehat A$. Finally, we write $$\Phi:A\to \widehat A, \,\,\, \widehat \Phi:\widehat A\to A$$ for the morphisms induced by $\Theta$ and $\widehat \Theta$, so that $$\Phi\circ \widehat \Phi=-1, \,\, \widehat \Phi\circ \Phi=-1.$$ 

Consider two orthogonal Mukai vectors $$v=(r, k\Theta, \chi), \,\, w=(r', k'\Theta, \chi'),$$ satisfying assumptions ${\bf (A.1)}$ and ${\bf (A.2)}.$ Central for our arguments is the following diagram \cite{Y}, \cite{abelian}: \begin{equation}\label{sdia}\xymatrix{K_{v}\times A \times \widehat
A\ar[r]^{\,\,\,\,\,\,\tau} \ar[d]^{p} & \uv \ar[d]^{\alpha} \\ A\times \widehat A \ar[r]^{\Psi} & A\times \widehat A}. \end{equation} Here, the morphism $$\tau:K_{v}\times A \times \widehat A\to \uv$$ is given by $$\tau(E, x, y)=t_{x}^{\star}E \otimes y.$$ Throughout this paper, we normalize the Albanese map $$\alpha=(\alpha^+, \alpha^-):\mathfrak M_v\to A\times \widehat A$$ of the introduction by \begin{equation}\label{normalb}\alpha(E)=\left(\det \fm ({E})\otimes \widehat \Theta^{-k}, \det E\otimes \Theta^{-k}\right).\end{equation} As before, we write $K_v$ for the fiber of $\alpha$ over the origin. Lemma $4.3$ of \cite {Y} identifies the lower horizontal morphism $$\Psi(x, y)=(-\chi x+k\widehat \Phi(y), k\Phi(x)+ry).$$ It can be seen that $\Psi$ is \'etale of degree $d_v^4=(k^2-r\chi)^4$, but we will not need this fact. 

\subsection {Properties of the theta bundles}\label{ptb}
As already mentioned, we prove the following analogue of the Dr\'ezet-Narashimhan theorem \cite{DN}, originally conjectured in \cite {abelian}. The reader will notice that the argument below applies in higher generality (in particular, in Section \ref{s3} we will need this result for nonprincipal polarizations).

\addtocounter{theorem}{-4}
\begin{theorem}\label{ct1} If $F_1, F_2$ have the same Mukai vector orthogonal to $v$, then $$\Theta_{F_1}=\Theta_{F_2}\otimes \left((-1)\circ \alpha^{+}\right)^{\star} \left(\det F_1\otimes \left(\det F_2\right)^{\vee}\right)\otimes \left((-1)\circ \alpha^{-}\right)^{\star}\left(\det \fm (F_1)\otimes \det \fm (F_2)^{\vee}\right).$$ 
\end{theorem}

\proof We begin by noting that $\Theta_F$ depends {\it a priori} on the holomorphic $K$-theory class of $F$. To compare the different theta bundles, we consider the virtual difference $$\mathsf f=F_1-F_2.$$ By Lemma $1$ of \cite {abelian} or alternatively by Lemma $17$ of \cite {BS}, in $K$-theory we can write \begin{equation}\label{bs}\mathsf f=M_1-M_2+\mathcal O_{Z_1}-\mathcal O_{Z_2},\end{equation} for line bundles $M_1, M_2$ over $A$, and zero dimensional subschemes $Z_1, Z_2$. (The proof in \cite {abelian} is written for sheaves, but the case of complexes is a consequence.) By assumption, the Mukai vector of $\mathsf f$ equals $0$. In particular $$c_1(M_1)=c_1(M_2)\implies \chi(M_1)=\chi(M_2).$$ Since $\chi(\mathsf f)=0$, it follows that the lengths of $Z_1$ and $Z_2$ are equal $$\ell(Z_1)=\ell(Z_2):=\ell.$$ 

Let us write $$\mathsf m=M_1-M_2, \,\,\,\mathsf n=\mathcal O_{Z_1}-\mathcal O_{Z_2},$$ and observe that $$\Theta_{\mathsf f}=\Theta_{\mathsf m}\otimes \Theta_{\mathsf n}.$$ To prove the Theorem, we show $$\Theta_{\mathsf m}\cong ((-1)\circ \alpha^+)^{\star} \left(\det \mathsf m\right) \boxtimes ((-1)\circ \alpha^-)^{\star} \left(\det \fm (\mathsf m)\right)$$ and $$\Theta_{\mathsf n}\cong\left((-1)\circ \alpha^{-}\right)^{\star}\det \fm ({\mathsf n}).$$

Note that $$\det \fm (\mathsf n)={\mathcal P}_{a(Z_1)-a(Z_2)},$$ where $a$ is the addition morphism. The isomorphism $$\Theta_{\mathsf n}\cong\left(\alpha^{-}\right)^{\star}\mathcal P_{a(Z_2)-a(Z_1)}$$ will be checked along any test family of sheaves. Indeed, consider a flat family $$\mathcal E\to S\times A$$ of sheaves of type $v$, inducing a morphism  $$\alpha^{-}:S\to \widehat A$$ by taking determinants and twisting. Remark that by the see-saw theorem 
 $$\det \mathcal E=(\alpha^{-}\times 1)^{\star} \mathcal P\otimes p^{\star} \mathcal V\otimes q^{\star} D,$$ for some line bundle $\mathcal V\to S$, with $p, q$ being the projections from $S\times A$ to $S$ and $A$. Here, we wrote $D$ for the line bundle used to twist the determinants to reach degree zero (of course, in our case, $D=\Theta^k$, but the notation emphasizes that the proof is general). For a zero dimensional subscheme $Z$ of length $\ell$, we calculate $$\det \rp p_{!} (\mathcal E\otimes q^{\star} \mathcal O_Z)=\bigotimes_{z\in Z} \det \mathcal E_{z}=\bigotimes_{z\in Z} \left((\alpha^{-})^{\star}\mathcal P_z\otimes \mathcal V\right)=(\alpha^{-})^{\star}\mathcal P_{a(Z)}\otimes  \mathcal V^{\ell}.$$ Therefore, $$\Theta_{\mathsf n}=\left((\alpha^{-})^{\star}\mathcal P_{a(Z_1)}\otimes \mathcal V^{\ell}\right)^{-1}\otimes \left((\alpha^{-})^{\star} \mathcal P_{a(Z_2)}\otimes \mathcal V^{\ell}\right)=\left(\alpha^{-}\right)^{\star}\mathcal P_{a(Z_2)-a(Z_1)},$$ as claimed. 

It remains to prove the first isomorphism.  To simplify notation, let us write $$M_1=M\otimes \zeta,\,\,\, M_2=M$$ for $\zeta\in \widehat A$ a degree $0$ line bundle, and $M$ an arbitrary line bundle over $A$. Set $$\widehat M=\det \fm (M).$$ We have $\det \mathsf m=\zeta,$ and $$\det \fm (\mathsf m)=\det \fm (M\otimes \zeta)\otimes \det \fm (M)^{-1}=t_{\zeta}^\star\widehat M\otimes \widehat M^{-1}=\widehat \Phi_{M}(\zeta)$$ where $$\widehat \Phi_{M}:\widehat A\to A$$ is the homomorphism induced by $\widehat M$. We show \begin{equation}\label{eqlty}\Theta_{\mathsf m}=(\alpha^{+})^{\star} \zeta^{-1}\boxtimes (\alpha^-)^{\star} \Phi_{\widehat M}(\zeta)^{-1}.\end{equation} First, we will verify this equality after pullback by $\tau$. To begin, we find the pullback of the right hand side: \begin{eqnarray*}\tau^{\star} \left((\alpha^+)^{\star} \zeta^{-1}\boxtimes (\alpha^-)^{\star} \Phi_{\widehat M}(\zeta)^{-1}\right)&=&\text{pr}^{\star}\Psi^{\star}\left(\zeta^{-1}\boxtimes \Phi_{\widehat M}(\zeta)^{-1}\right)\\&=&\text{pr}^{\star} \left(-\chi x + k\widehat \Phi(y), k\Phi(x)+ry\right)^{\star} \left(\mathcal \zeta^{-1}\boxtimes \Phi_{\widehat M}(\zeta)^{-1}\right)\\&=&\text{pr}^{\star}\left(\zeta^{\chi}\otimes \Phi^{\star}(\widehat \Phi_{M}(\zeta))^{-k}\boxtimes \widehat \Phi^{\star}{\zeta}^{-k}\otimes \widehat \Phi_M(\zeta)^{-r} \right).\end{eqnarray*} We claim next that $$\tau^{\star}\Theta_{\mathsf m}=\mathcal Q\boxtimes \text{pr}^{\star}\left(\zeta^{\chi}\otimes \Phi^{\star}(\widehat \Phi_{M}(\zeta))^{-k}\boxtimes \widehat \Phi^{\star}{\zeta}^{-k}\otimes \widehat \Phi_M(\zeta)^{-r} \right)$$ for some line bundle $\mathcal Q$ over $K_v$, where $\mathcal Q$ is the restriction of $\tau^{\star} \Theta_{\mathsf m}$ to $K_v$. For $M$ fixed, the Chern class of $\mathcal Q$ depends only on the Chern class of $\zeta$, and since for $\zeta=\mathcal O$ we obtain the trivial bundle, this must be the case for all $\mathcal \zeta$'s of degree zero. Since $K_v$ is simply connected, $\mathcal Q$ must be trivial. 

To prove the splitting above, note that for $E\in K_v$, the restriction of $\tau^{\star}\Theta_{\mathsf m}$ to $\{E\}\times A\times \widehat A$ equals $$\mathcal L_E=\det {\mathbf R} p_{23!}(m_{12}^{\star} E\otimes p_{13}^{\star}\mathcal P\otimes p_1^{\star} (M\otimes \zeta-M))^{-1},$$ where $$m_{12}:A\times A\times \widehat A\to A$$ denotes the addition on the first two factors, and the $p$'s are the projections over the factors of $A\times A\times \widehat A$. We prove that \begin{equation}\label{ele}\mathcal L_E=\left(\chi(E)\cdot \zeta-\Phi_{D}(\widehat \Phi_{M}(\zeta))\right)\boxtimes \left(-\widehat \Phi_D ({\zeta})-\text{rank }(E)\cdot \widehat \Phi_{M}(\zeta)\right),\end{equation} where $D$ is the determinant of $E$ (of course, $D=\Theta^k$, but as before we prefer to write a general proof of the Theorem). Here $$\Phi_D:A\to \widehat A,\,\,\, \widehat \Phi_D:\widehat A\to A$$ are the Mumford homomorphisms induced by the line bundles $D$ and $\det \fm (D)$ respectively. We also switched to additive notation, for ease of reading. 

The idea of the proof is already contained in the above argument. 
We first note that $\mathcal L_E$ only depends on the holomorphic $K$-theory class of $E$. In fact, we argue that $\mathcal L_E$ depends on the rank, determinant and Euler characteristic of $E$.  Indeed, for two sheaves $E_1, E_2$ with the same data as above, we form the virtual difference $$\mathsf e=E_1-E_2.$$ Using equation \eqref{bs}, we can write $$\mathsf e=\mathcal O_Z-\mathcal O_W,\,\,\, \text { with } \ell (Z)=\ell(W).$$ Then $$\mathcal L_{\mathsf e}=\mathcal L_{E_1}\otimes \mathcal L_{E_2}^{-1}$$ is trivial, since for any $a\in A$, we have \begin {eqnarray*}\det {\mathbf R} p_{23!} (m_{12}^{\star} \mathcal O_a&\otimes& p_{13}^{\star} \mathcal P\otimes p_{1}^{\star} (M\otimes \zeta-M))\\&=&\det \left((t_{a}\times 1)\circ (-1 \times 1)\right)^{\star}(\mathcal P\otimes p_{1}^{\star} (M\otimes \zeta-M))\\&=&\left((t_{a}\times 1)\circ (-1 \times 1)\right)^{\star}\,\, p_1^{\star}\zeta=\text{pr}^{\star} \mathcal \zeta^{-1}.\end{eqnarray*} 
Since for fixed rank and Euler characteristic, $\mathcal L_E$ depends only on the determinant of $E$, we may assume that $E$ splits into rank $1$ and rank $0$ factors. Since both sides of \eqref{ele} change multiplicatively as $E$ splits, it therefore suffices to consider the cases $$E=\text{line bundle},\,\, \,\, E=\mathcal O_{a},\,\, a\in A.$$ 

Now, note that for $E=\mathcal O_a$ we obtain by the preceding paragraph that $$\mathcal L_{\mathcal O_a}=\det {\mathbf R} p_{23!}(m_{12}^{\star} \mathcal O_a\otimes p_{13}^{\star} \mathcal P\otimes p_1^{\star} (M\otimes \zeta-M))^{-1}\cong \text{pr}_A^{\star}{\zeta}.$$ The calculation when $E=\mathcal O(D)$ is a line bundle will be more involved. We need to show $$\mathcal L=\left(\chi(D)\cdot \zeta-\Phi_{D}(\widehat \Phi_{M}(\zeta))\right)\boxtimes \left(-\widehat \Phi_D({\zeta})-\widehat \Phi_{M}(\zeta)\right).$$ 
Observe that over $A\times A$ we have $$m_{12}^{\star} \mathcal O(D)=(1\times \Phi_D)^{\star} \mathcal P\otimes p_{1}^{\star} \mathcal O(D)\otimes p_2^{\star} \mathcal O(D).$$ We compute \begin{eqnarray*}\mathcal L^{-1}&=&\det {\mathbf R}p_{23!\,} (p_{12}^{\star}(1\times \Phi_D)^{\star} \mathcal P\otimes p_{13}^{\star} \mathcal P \otimes p_1^{\star} (M\otimes \zeta\otimes D-M\otimes D))\\&=&(\Phi_D \times 1)^{\star} \det {\mathbf R}p_{23!\,} (p_{12}^{\star} \mathcal P\otimes p_{13}^{\star} \mathcal P\otimes p_1^{\star} (M\otimes \zeta\otimes D-M\otimes D)).\end{eqnarray*} The pullback $p_2^{\star} \mathcal O(D)$ did not contribute above by the projection formula and the vanishing of the Euler characteristic. Note furthermore that $$p_{12}^{\star} \mathcal P\otimes p_{13}^{\star} \mathcal P=(1\times \widehat m)^{\star} \mathcal P$$ where $$1\times \widehat m: A\times \widehat A\times \widehat A\to A\times \widehat A$$ is the addition map. We conclude \begin{eqnarray*}\mathcal L^{-1}&=&(\Phi_D \times 1)^{\star} \det {\mathbf R}p_{23!\,} \left(\left(1\times \widehat m\right)^{\star} \mathcal P\otimes p_1^{\star} \left(M\otimes D\otimes \zeta-M\otimes D\right)\right)\\&=&(\Phi_D \times 1)^{\star} \widehat m^{\star}\, \det {\mathbf R}p_{2!\,} \left(\mathcal P\otimes p_{1}^{\star} (M\otimes D\otimes \zeta-M\otimes D)\right) \\ &=& (\Phi_D \times 1)^{\star} \widehat m^{\star} \left(\det \fm ({M\otimes D\otimes \zeta})\otimes \det \fm (M\otimes D)^{-1}\right)\\&=&(\Phi_D \times 1)^{\star} \widehat m^{\star} \,\left(t_{\mathcal \zeta}^{\star} \det \fm (M\otimes D)\otimes \det \fm (M\otimes D)^{-1}\right)\\&=&(\Phi_D \times 1)^{\star} \widehat m^{\star} \,{\widehat \Phi_{M+D}} (\zeta).\end{eqnarray*} Now for a degree $0$ line bundle $U$ we have $$\widehat m^{\star} U=p_1^{\star} U\otimes p_{2}^{\star} U,$$ hence \begin{eqnarray*}\mathcal L^{-1}&=&(\Phi_D \times 1)^{\star}\left(\widehat{\Phi}_{D+M}\,  (\zeta)\boxtimes \widehat {\Phi}_{D+M}\, (\zeta)\right)\\&=&(\Phi_D \times 1)^{\star}\left((\widehat \Phi_{D}\,(\zeta)+\widehat {\Phi}_M (\zeta))\boxtimes (\widehat \Phi_{D}\,(\zeta)+\widehat {\Phi}_M (\zeta))\right)\\&=&\left(-\chi(D)\cdot \zeta +\Phi_{D}(\widehat{\Phi}_{M}(\zeta))\right) \boxtimes \left(\widehat{\Phi}_D (\zeta)+\widehat {\Phi}_M(\zeta)\right)\end{eqnarray*} as claimed in \eqref{ele}. For the last equality, we used that $\Phi_{D}\circ \widehat \Phi_D=-\chi(D) \cdot {\mathbf 1}$, see Lemma $4.2$ of \cite {Y}.

Equality \eqref{eqlty} is now checked under pullback by $\tau$. To complete the argument, fix $M$. Observe that the assignment $$\zeta\to \Theta_{\mathsf m}\otimes \left((\alpha^{+})^{\star} \zeta\boxtimes (\alpha^-)^{\star} \Phi_{\widehat M}(\zeta)\right)$$ defines a morphism $$\pi:\widehat A\to \text{Pic } (\uv).$$ Note moreover that the above discussion implies that $$\tau^{\star}\circ \pi=0.$$Since the kernel of $\tau^{\star}$ is discrete, $\pi$ must be constant. Since $\pi(\mathcal O)=\mathcal O$, we must have $\pi(\zeta)=\mathcal O$ for all $\zeta\in \widehat A$, completing the proof. \qed
\vskip.1in

\begin {convention}\label{r1} The theorem above shows that $\Theta_F$ only depends on the rank, Euler characteristic, determinant and determinant of the Fourier-Mukai of the bundle $F$. The Mukai vectors used in this paper are all of the form $$w=(r', k'\Theta, \chi').$$ We define the normalized theta bundles $$\Theta_w:=\Theta_F\to \uv,$$ for complexes $F$ satisfying \begin{equation}\label{conv1}\text {rank } F=r', \,\, \chi(F)=\chi',\,\, \det F=\Theta^{k'}, \,\, \det \fm (F)=\widehat \Theta^{k'}.\end{equation} Even though the exact choice for $F$ is irrelevant, for concreteness we may take $$F=(r'-1)\mathcal O\oplus \Theta^{k'}\oplus (\chi'-k'^2)\mathcal O_{o},$$ with $o\in A$ denoting the origin. The normalization we use is aligned with that of the Albanese morphism in \eqref{normalb}.  \end{convention}

\begin {example}\label{example1} Assume that $v=(1, 0, -n)$ so that $\uv\cong A^{[n]}\times \widehat A$ via the isomorphism $$(Z, y)\mapsto I_Z\otimes y.$$ Then, $\alpha$ becomes the morphism $$(-a, 1):A^{[n]}\times \widehat A\to A\times \widehat A,$$ where as usual $a$ is the addition map. For a sheaf $F\to A$ of rank $r$, we obtain $$\Theta_F=\det {\mathbf R}p_{12!} (p_{13}^{\star}I_{\mathcal Z}\otimes p_3^{\star}F\otimes p_{23}^{\star}\mathcal P)^{-1},$$ where the projections are considered onto the factors of the product $A^{[n]}\times \widehat A\times A$ and $\mathcal Z$ is the universal subscheme in $A^{[n]}\times A$. This yields $$\Theta_F=\det \mathbf Rp_{12!}(p_3^{\star} F\otimes p_{23}^{\star} \mathcal P)^{-1}\otimes \det \mathbf Rp_{12!} (p_{13}^{\star}\mathcal O_{\mathcal Z} \otimes p_3^{\star}F \otimes p_{23}^{\star} \mathcal P).$$ The second line bundle can be found via the see-saw theorem and Section $5$ of \cite{EGL} $$(\det F)_{(n)}\otimes M^{r}\otimes (a,1)^{\star} \mathcal P^r.$$ Here, $M$ is half the exceptional divisor on the Hilbert scheme and $(\cdot)_{(n)}$ denotes the symmetrization of a line bundle from $A$ over $\text{Sym}^n(A)$. Therefore, over $A^{[n]}\times \widehat A$, $$\Theta_F=\left(\left((\det F)_{(n)}\otimes M^r\right)\boxtimes \det \fm(F)^{-1}\right)\otimes (a,1)^{\star}\mathcal P^{r}.$$ This expression is consistent with the statement of the theorem. 

\end {example}

\subsection {Theta bundles and \'etale pullbacks} \label {plbckk} We now return to the \'etale diagram of Section \ref{setup}: \begin{center} $\xymatrix{K_{v}\times A \times \widehat
A\ar[r]^{\,\,\,\,\,\,\tau} \ar[d]^{p} & \uv \ar[d]^{\alpha} \\ A\times \widehat A \ar[r]^{\Psi} & A\times \widehat A}.$ \end{center} We continue to refer the reader to \cite {abelian}, where we showed the splitting of the pullback $$\tau^{\star} \Theta_w=\Theta_w\boxtimes \mathcal L.$$ The line bundle $\Theta_w$ was shown to be independent of choices on the simply connected manifold $K_v$. Furhermore, we proved that $$\mathcal L=\left ( \det \rp p_{23!}\left(m_{12}^{\star}E\otimes 
p_{13}^{\star}\mathcal
P\otimes p_{1}^{\star}F\right) \right)^{-1},$$ for $E\in K_v$ and  for $F$ satisfying \eqref{conv1}. Here, the $p$'s denote the projections on the corresponding factors of $A\times A\times \widehat A,$ while 
$$m_{12}:A\times
A\times \widehat A\to A$$ is as usual the addition on the first two factors.  We moreover calculated the Euler characteristic of $\mathcal L$ in \cite {abelian}. In the lemma below, we identify $\mathcal L$ explicitly.  As a corollary we obtain \begin{equation}\label{pback}\Psi^{\star} \mathbf E(v, w)=p_{\star} \tau^{\star} \Theta_w=H^0(K_v, \Theta_w)\otimes \mathcal L.\end{equation} 
 
\begin {lemma} \label{lemmal}We have $$\mathcal L=\left(\Theta^{-\chi'k-\chi k'}\boxtimes \widehat\Theta^{-rk'-r'k} \right)\otimes \mathcal P^{r'\chi+kk'}.$$ \end{lemma}

The Lemma will be used in Section \ref{vbb} for $k=0$, but for future reference, we write the proof for all $k$. The result is also valid for nonprincipal polarizations. In this case, the exponent of the Poincar\'e bundle should be modified to $r'\chi+kk'e$ where $e=\chi(\Theta)$. 

\proof  The proof is similar to that of Theorem \ref{tpt}. 
Pick two complexes $E$ and $F$ satisfying Convention $1$, not necessarily orthogonal in the Mukai pairing. First, we note that the bundle $$\mathcal L_{E, F}=\left ( \det \rp p_{23!}\left(m_{12}^{\star}E\otimes 
p_{13}^{\star}\mathcal
P\otimes p_{1}^{\star}F\right) \right)^{-1}$$ depends only on the holomorphic $K$-theory classes of $E$ and $F$. We will furthermore remark below that the line bundle only depends on the Mukai vectors $v$ and $w$, the determinant and determinant of the Fourier-Mukai of $E$ and $F$. In fact, only the statement about the first argument $E$ will be useful to us, so that we show $$\mathcal L_{E, F}\cong \mathcal L_{E^{'}, F}$$ when $(E, E^{'})$ have the same Mukai vectors, determinants and determinants of Fourier-Mukai. 

To prove this isomorphism, consider the virtual sheaf $$\mathsf e=E-E^{'}.$$ Note that by equation \eqref{bs}, in $K$-theory we have $$\mathsf e=\mathcal O_Z-\mathcal O_W$$ for two zero-dimensional subschemes which have the same length. Since the Fourier-Mukai transform of $\mathsf e$ has trivial determinant, we must have $a(Z)=a(W).$ We prove $$\mathcal L_{\mathsf e, F}=\mathcal L_{E, F}\otimes \mathcal L_{E^{'}, F}^{-1}$$ is trivial. Over $A\times \widehat A$, we calculate \begin{eqnarray*}\det \rp p_{23!} (m_{12}\mathcal O_Z\otimes p_{13}^{\star}\mathcal P\otimes p_1^{\star} F)&\cong& \bigotimes_{z\in Z} ((t_{z}\times 1)\circ (-1\times 1))^{\star}\det \left (\mathcal P\otimes p_1^{\star} F\right)\end{eqnarray*} which only depends on $a(Z)$ by the theorem of the square. Thus, we obtain the same answer replacing $Z$ by $W$, therefore showing $\mathcal L_{\mathsf e, F}$ is trivial. 

With this understood, we prove the lemma. We may assume then that $E$ splits as a direct sum of copies of $\mathcal O$, $\Theta$ and structure sheaves $\mathcal O_{o}$ $$E=(r-k)\mathcal O+k\Theta+(\chi-k)\mathcal O_o.$$ In fact, since $\mathcal L_{E, F}$ is multiplicative as $E$ splits, it suffices to prove the lemma separately for the three sheaves $$E=\mathcal O, \,\, E= \Theta\,\, \text{ and } E=\mathcal O_{o}.$$ First, for $E=\mathcal O$, we obtain $$\mathcal L=\det \rp p_{23!} (m_{12}^{\star} \mathcal O\otimes p_{13}^{\star}\mathcal P\otimes p_1^{\star} F)^{-1}\cong \mathcal O\boxtimes \left(\det \fm (F)\right)^{-1}=\mathcal O\boxtimes \widehat \Theta^{-k'},$$ while for $E=\mathcal O_o$, we have $$\mathcal L=\det \rp p_{23!} (m_{12} \mathcal O_{o}\otimes p_{13}^{\star}\mathcal P\otimes p_1^{\star} F)^{-1}\cong \det \,(-1, 1)^{\star} \left(\mathcal P\otimes p_1^{\star} F\right)^{-1}\cong \mathcal P^{r'}\otimes \Theta^{-k'}.$$ The calculation for $E=\Theta$ is more involved.  We show $$\mathcal L=\det \rp p_{23!} (m_{12}^{\star} \Theta \otimes p_{13}^{\star} \mathcal P\otimes p_1^{\star} F)^{-1}=\left(\Theta^{-\chi'-k'}\boxtimes \widehat \Theta^{-k'-r'}\right)\otimes \mathcal P^{r'+k'}.$$ Observe that over $A\times A$ $$m^{\star} \Theta=(1\times \Phi)^{\star} \mathcal P\otimes p_1^{\star} \Theta\otimes p_2^{\star} \Theta.$$
We calculate \begin{eqnarray*}\mathcal L&=&\det \rp p_{23!}((1\times \Phi)^{\star}p_{12}^{\star}\mathcal P\otimes p_{13}^{\star} \mathcal P\otimes p_1^{\star} (F\otimes \Theta)\otimes p_2^{\star}\Theta)^{-1}\\&=&\det \rp p_{23!}((1\times \Phi)^{\star}p_{12}^{\star}\mathcal P\otimes p_{13}^{\star} \mathcal P\otimes p_1^{\star} (F\otimes \Theta))^{-1}\otimes \text{pr}_1^{\star}\,\Theta^{-r'-\chi'-2k'},\end{eqnarray*} where we noted that $$\chi(F\otimes \Theta)=r'+\chi'+2k'.$$ This in turn becomes 
\begin{eqnarray*} \mathcal L
&=&(\Phi\times 1)^{\star}\det \rp p_{23!} (p_{12}^{\star}\mathcal P\otimes p_{13}^{\star} \mathcal P\otimes p_1^{\star} (F\otimes \Theta))^{-1}\otimes \text{pr}_1^{\star}\Theta^{-r'-\chi'-2k'}\\&=&(\Phi\times 1)^{\star} \det \rp p_{23!} ((1\times \widehat m)^{\star} \mathcal P\otimes p_1^{\star} (F\otimes \Theta))^{-1}\otimes \text{pr}_1^{\star}\,\Theta^{-r'-\chi'-2k'}\end{eqnarray*} where we noted again that over $A\times \widehat A\times \widehat A$ we have $$p_{12}^{\star}\mathcal P\otimes p_{13}^{\star}\mathcal P=(1\times \widehat m)^{\star} \mathcal P.$$ 
We continue the calculation 
\begin{eqnarray*}\mathcal L&=&(\Phi\times 1)^{\star} \det \left(\rp p_{23!}((1\times \widehat m)^{\star}\mathcal P\otimes p_1^{\star}(F\otimes \Theta))\right)^{-1}\otimes \Theta^{-r'-\chi'-2k'}\\ &=&(\Phi\times 1)^{\star} \widehat m^{\star} \det \rp p_{3!}\left(\mathcal P\otimes p_1^{\star} (F\otimes \Theta)\right)^{-1} \otimes \Theta^{-r'-\chi'-2k'}\\&=&(\Phi\times 1)^{\star}\widehat m^{\star} \det \fm (F\otimes \Theta)^{-1}\otimes \Theta^{-r'-\chi'-2k'}\\&=&(\Phi\times 1)^{\star} \widehat m^{\star} \widehat \Theta^{-r'-k'}\otimes \Theta^{-r'-\chi'-2k'}\end{eqnarray*} 
Noting again that $$\widehat m^{\star} \widehat \Theta= (\widehat \Phi\times 1)^{\star} \mathcal P\otimes p_1^{\star} \widehat \Theta\otimes p_2^{\star} \widehat \Theta,$$ we obtain the result \begin{eqnarray*}\mathcal L&=&(\Phi\times 1)^{\star} (\widehat \Phi\times 1)^{\star} \mathcal P^{-r'-k'}\otimes \left(\Phi^{\star} \widehat \Theta^{-r'-k'}\boxtimes \widehat \Theta^{-r'-k'}\right)\otimes \Theta^{-r'-\chi'-2k'}\\&=&(-1, 1)^{\star}\mathcal P^{-r'-k'}\otimes \left(\Theta^{r'+k'}\boxtimes \widehat \Theta^{-r'-k'}\right)\otimes \Theta^{-r'-\chi'-2k'}\\&=&\mathcal P^{r'+k'}\otimes \left(\Theta^{-\chi'-k'}\boxtimes \widehat \Theta^{-r'-k'}\right).\end{eqnarray*} The lemma is now proved. 

The only detail that still needs clarification is the fact, used above, that \begin{equation}\label{deteq}\det \fm ({F\otimes \Theta})=\widehat \Theta^{r'+k'}.\end{equation} Here, we crucially use that $\Theta$ is symmetric so that $$\det \fm ({\Theta^{k}})=\widehat \Theta^k.$$ This statement follows for instance by taking determinants in Lemma 2(ii) in \cite {wirtinger}. The lemma is stated for odd numerics, but the case $k$ even follows since both sides depend polynomially in $k$. For the left hand side, this is a general statement about integral transforms, which can be proved via induction on dimension. The inductive step consists in cutting with hyperplanes in the (pluri)-theta series to reduce dimension. 

To prove \eqref{deteq}, let $$\mathcal A_F=\det  \fm (F\otimes \Theta),$$ and observe that $\mathcal A_F$ depends on the rank, Euler characteristic, determinant and determinant of the Fourier-Mukai transform of $F$. Indeed, if $F_1$ and $F_2$ are two such sheaves, we write $$\mathsf f=F_1-F_2=\mathcal O_Z-\mathcal O_W,$$ where $a(Z)=a(W).$ But then $$\mathcal A_{F_1}\otimes \mathcal A_{F_2}^{-1}=\mathcal A_{\mathsf f}=\det \fm (\Theta\otimes \mathsf f)=\det \fm(\Theta \otimes (\mathcal O_Z-\mathcal O_W))=\mathcal P_{a(Z)-a(W)}=\mathcal O.$$
Therefore, it suffices to assume that $$F=(r'-1)\mathcal O\oplus \Theta^{k'}\oplus ({\chi'-k'^2})\mathcal O_{o}.$$ To conclude, note that in this case $$\mathcal A_{F}=\det \fm (\Theta)^{r'-1}\otimes \det \fm ({\Theta^{k'+1}})\otimes \det \fm(\Theta\otimes \mathcal O_o)^{\chi'-k'^2}=\widehat \Theta^{r'-1}\otimes \widehat \Theta^{k'+1}=\widehat \Theta^{r'+k'}.$$  \qed

\subsection {Construction of the strange duality map} \label{sdsection}In this subsection, we use Theorem \ref{tpt} to construct the duality map $\mathsf {SD}$  mentioned in equation \eqref{sdm} of the Introduction when $c_1(v\otimes w)\neq 0$. A similar construction was achieved in \cite {Po} in the case of curves by packaging together all the strange duality morphisms relatively over the Jacobian.

We assume that $$c_1(v\otimes w)\cdot \Theta>0.$$ The case $c_1(v\otimes w)\cdot \Theta<0$ is similar, and the needed modifications are explained in \cite{abelian}. Serre duality implies that for any two stable sheaves $E$ and $F$ we have $$H^2(E\otimes F)=0.$$ Furthermore, the locus \begin{equation}\label{thl}\mathcal D=\{(E, F): h^0(E\otimes F)\neq 0\}\hookrightarrow \uv\times \uw\end{equation} has expected codimension $1$. 
The defining equation of \eqref{thl} is used to prove that:

\begin {lemma} There exists a natural morphism $$\mathsf {SD}: \mathsf j^{\star}\mathbf E(v, w)^{\vee}\to \widehat {\mathbf E(w, v)}$$
\end {lemma}
Here we write $$\mathsf j:A\times \widehat A\to A\times \widehat A$$ for the multiplication by $(-1, -1)$. Since our Mukai vectors are in fact assumed to be symmetric, pullback by the morphism $\mathsf j$ will not be necessary. 

Note that both sides are locally free. For the left hand side, we use \eqref{pback}. For the right hand side, we argue that $\mathbf E(w, v)$ satisfies the index theorem with index $0$. This can be checked after pullback under isogenies. To this end, we invoke \eqref{pback} and the fact that under our assumptions, $\mathcal L$ is ample, cf. the criterion used in the proof of Theorem \ref{tttt}, {\it Step 3}.

\proof To construct $\mathsf {SD}$, we need a natural section of the bundle $$\mathsf  j^{\star}\mathbf E(v, w)\otimes \widehat {\mathbf E(w, v)}=\text{pr}_{12\star}(j^{\star}\mathbf E(v, w)\boxtimes \mathbf E(w, v)\otimes \mathcal Q)$$ where $$\text{pr}_{12}:A\times \widehat A\times A\times \widehat A\to A\times \widehat A$$ is the projection onto the first two factors, and $\mathcal Q\to (A\times \widehat A)\times (A\times \widehat A)$ is the Poincar\'e bundle on the self-dual abelian variety $A\times \widehat A$. Note furthermore that $$j^{\star}\mathbf E(v, w)\boxtimes \mathbf E(w, v)\otimes \mathcal Q=\alpha_{\star} (\Theta_w\boxtimes \Theta_v\otimes \alpha^{\star}\mathcal Q)$$ where $$\alpha =(\mathsf j\circ \alpha_v)\times \alpha_w: \uv\times \uw\to (A\times \widehat A)\times (A\times \widehat A).$$ We will construct a natural section of the line bundle $$\Theta_w\boxtimes \Theta_v\otimes \alpha^{\star}\mathcal Q\to \mathfrak M_v\times \mathfrak M_w,$$ which we will then pushforward by $\text{pr}_{12}\circ \alpha$ to complete the proof. 

For simplicity let us assume that universal sheaves $\mathcal E\to \uv\times A$ and $\mathcal F\to \uw\times A$ exist, using quasi-universal families otherwise. We form the bundle $$\mathcal D=\det \mathbf Rp_{!}(\mathcal E\boxtimes^{\mathbf L} \mathcal F)^{-1},$$ obtained by pushforward via the projection $$p:\uv\times \uw\times A\to \uv\times \uw.$$ 
We claim that $$\mathcal D=\Theta_w\boxtimes \Theta_v\otimes \alpha^{\star} \mathcal Q.$$ This is precisely the see-saw principle combined with Theorem \ref{ct1}. Indeed, the restriction $\Theta_E$ of $\Theta$ to $\{E\}\times \uw$ equals $$\Theta_E=\Theta_{v}\otimes \alpha_w^{\star}\left( (\det E\otimes \Theta^{-k})^{\vee}\boxtimes (\det {\fm (E)}\otimes \widehat \Theta^{-k})^{\vee}\right)=\Theta_v\otimes \alpha^{\star} \mathcal Q|_{\{E\}\times \uw}.$$ The calculation of the restriction to $\uv\times \{F\}$ is similar: $$\Theta_{F}=\Theta_w\otimes \alpha_v^{\star}\left((\det F\otimes \Theta^{-k'})^{\vee}\boxtimes (\det \fm (F)\otimes \widehat \Theta^{-k'})^{\vee}\right)=\Theta_w\otimes \alpha^{\star} \mathcal Q|_{\uv\times \{F\}}.$$

We can now complete the proof. The locus where $$Tor^1(E, F)=Tor^2(E, F)=0$$ has complement of codimension at least $2$ in the product space by Proposition $0.5$ in \cite {Y2}; see also Remark $4.4$ in \cite {Y}.  Indeed, the Proposition guarantees that when both ranks are at least $2$, $E$ and $F$ are generically locally free in their corresponding moduli spaces, giving the claim. The only exception occurs for vectors of the form $$v=r\exp(c_1(L))-[\text{pt}]$$ for a fixed line bundle $L\to A$. In this case, $E$ is always nonlocally free. In fact, by Corollary $4.4$ in \cite {mukai3}, all such $E$'s  are obtained as kernels of morphisms $$H\otimes L\to \mathbb C_x$$ for $x\in A$, and for homogeneous vector bundles $H\to A$, so $E$ fails to be locally free at $x$. The codimension analysis holds true in this case as well: for each fixed $F$, the locus of $E$'s which intersects the singularity locus of $F$ has codimension $2$ in the moduli space. The argument when either $E$ or $F$ have rank $1$ also follows by a moduli count.  

Now, along the locus of vanishing Tor's, the pushforward $\mathbf Rp_{!}(\mathcal E\boxtimes^{\mathbf L} \mathcal F)$ can be represented by a two step complex of equal rank vector bundles $$0\to \mathcal A_0\stackrel{\sigma}{\to} \mathcal A_1\to 0.$$ This yields a section $\det \sigma$ of $$\det \mathcal A_1\otimes \det \mathcal A_0^{-1}=\det \mathbf Rp_{!}(\mathcal E\boxtimes^{\mathbf L} \mathcal F)^{-1}=\mathcal D=\Theta_w\boxtimes \Theta_v\otimes \alpha^{\star}\mathcal Q$$ vanishing precisely along the theta locus \eqref{thl}. \qed

\begin{conjecture} \label{cj1}The morphism $$\mathsf {SD}:\mathsf j^{\star}\mathbf E(v,w)^{\vee}\to \widehat {\mathbf E(w,v)}$$ is an isomorphism. 
\end{conjecture}

\begin {remark} Just as in the case of curves, there is a slight asymmetry in the roles of $v$ and $w$ in the strange duality morphism \eqref{sdi}: on one side, the determinant and determinant of the Fourier-Mukai vary, while on the other side these invariants are kept fixed. However, just as in the case of curves \cite {Pol}, the above reformulation makes it clear that
\begin {corollary} If the duality morphism \eqref{sdi} is an isomorphism for the pair $(v, w)$, then it is an isomorphism for the pair $(w, v)$.
\end {corollary}

\end {remark}

\section {The Verlinde bundles in degree $0$} \label{s3}

This section contains the proof of the main results. 
The crux of the matter is the Lefschetz-Riemann-Roch trace calculation of Theorem \ref{t1}. This in turn implies Theorem \ref{thm21}. We also prove the Fourier-Mukai symmetries of the Verlinde bundles in Theorem \ref{tttt}, thus offering evidence for Conjecture \ref{cj1}. 

\subsection {Setup} We specialize to the case $c_1(v)=0$ {\it i.e.} we assume $$v=(r, 0, \chi), \,\, w=(rh, k\Theta, -\chi h)$$ with $r, \chi$ odd, and $(r, \chi)=1$. The integers $(rh, k, -\chi h)$ were previously denoted $(r', k', \chi')$, but the new notation should make the exposition easier to follow. 

For these numerics, the usual \'etale diagram \begin{center} $\xymatrix{K_{v}\times A \times \widehat
A\ar[r]^{\,\,\,\,\tau} \ar[d]^{p} & \uv \ar[d]^{\alpha} \\ A\times \widehat A \ar[r]^{\Psi} & A\times \widehat A}.$ \end{center} takes a simpler form. In particular, $$\Psi(x, y)=(-\chi x, ry).$$ In Subsection \ref{plbckk}, we calculated the pullback \begin{equation}\label{pback2}\tau^{\star} \Theta_w=\Theta_w \boxtimes \left(\Theta^{-\chi}\boxtimes \widehat\Theta^{-r} \right)^{k}\otimes \left(\mathcal P^{r\chi}\right)^{ h}=\Theta_w \boxtimes \left(\Theta^{-\chi}\boxtimes \widehat\Theta^{-r} \right)^{k}\otimes \Psi^{\star} \mathcal P^{ -h}.\end{equation}
 Consequently, we have \begin{equation}\label{pback1}\Psi^{\star} \mathbf E(v, w)=p_{\star} \tau^{\star} \Theta_w=H^0(K_v, \Theta_w)\otimes \left(\Theta^{-\chi}\boxtimes \widehat\Theta^{-r} \right)^{k}\otimes \Psi^{\star} \mathcal P^{-h}.\end{equation} 

\subsection {Group actions}\label{ga} We consider actions of the torsion group $$\mathsf G=A[-\chi]\times \widehat A[r]$$ on the three spaces $K_v$, $A\times \widehat A$ and $\uv$ appearing in the diagram. The action of $\mathsf G$ on $A\times \widehat A$ is given by translation on both factors, the action on $\uv$ is trivial, while the action on $K_v$ is given by $$(x, y)\cdot E=t_{-x}^{\star} E\otimes y^{-1}.$$ 

There are induced actions of a certain subgroup $\mathsf K$ of $\mathsf G$ on the theta bundles, which we now describe. Writing $$a=\gcd(k, \chi),\,\, b=\gcd(k, r),$$ we conclude that $$\gcd (a,b)=1 \text { and } a, b \text { are odd}.$$ We set $$\mathsf K=A[a]\times \widehat A[b]\hookrightarrow \mathsf G.$$ The morphisms $\tau$ and $\Psi$ are invariant under the action of $\mathsf G$, hence also under the action of $\mathsf K$. The bundles $\tau^{\star}\Theta_w$ and $\Psi^{\star}\mathcal P$ in \eqref{pback2} are naturally $\mathsf K$-equivariant.  We claim that $$\Theta_w\to K_v$$ is $\mathsf K$-equivariant as well.

Certainly, $\Theta^{-\chi k}\boxtimes \widehat \Theta^{-rk}$ carries a linearization of the  Heisenberg group $$\mathsf H=\mathsf H[-\chi k]\times \widehat {\mathsf H}[r k]$$ defined as an extension $$0\to \mathbb C^{\star} \times \mathbb C^{\star} \to \mathsf H\to A[-\chi k]\times \widehat A[rk]\to 0.$$ Our convention is that $\mathsf H[m]$ denotes the Heisenberg group of $\Theta^m$ consisting of pairs $(x, f)$ where $$f:t^{\star}_x \Theta^m\to \Theta^m$$ is an isomorphism, while $\widehat {\mathsf H}[n]$ denotes the Heisenberg group of $\widehat \Theta^{n}$. We have a natural morphism $$\iota:\mathsf H[a]\times \widehat {\mathsf H}[b]\to \mathsf H[-\chi k]\times \widehat {\mathsf H}[r k]=\mathsf H$$ which over the center restricts to $$(\alpha, \beta)\to (\alpha^{-\frac{\chi k}{a}}, \beta^{\frac{r k}{b}}).$$ It is now useful to pass to the finite Heisenberg groups $\widetilde {\mathsf H}$ obtained by restricting the centers to roots of unity. For instance $$0\to \mu_a\to \widetilde {\mathsf H}[a]\to A[a]\to 0,\,\, 0\to \mu_b\to \widetilde{\mathsf H}[b]\to \widehat A[b]\to 0.$$  Since the center of $\mathsf {\widetilde H}[a]\times \mathsf {\widetilde H}[b]$ is trivial under $\iota$, we obtain a morphism $$j:A[a]\times \widehat A[b]\to {\mathsf H}.$$ Furthermore, since $a, b$ are odd, the identification $$\Theta^{-\chi k}\boxtimes \widehat \Theta^{-rk}\cong (a,b)^{\star} \left(\Theta^{\frac{-\chi k}{a^2}}\boxtimes \widehat \Theta^{-\frac{k\chi}{b^2}}\right)$$ is compatible with the action of the group $$\mathsf K=A[a]\times \widehat A[b]$$ coming from $j$ on the left, and via the pullback on the right.  Using \eqref{pback2}, and the natural $\mathsf K$-action on line bundles $\tau^{\star} \Theta_w$ and $\Psi^{\star} \mathcal P$, we obtain a $\mathsf K$-linearization of $$\Theta_w\to K_v.$$ A similar construction was carried out in Section $3.1.2$ of \cite {wirtinger}. We will see below that the linearization of $\Theta_w\to K_v$ has the property that each $\zeta\in \mathsf K$ acts trivially in the fibers over the $\zeta$-fixed points.
  
\subsection{Trace calculations} \label{trmain} As a consequence of the above discussion, the signed sums of cohomologies $$\chi(K_v, \Theta_w)=\sum_i (-1)^i H^i(K_v, \Theta_w)$$ carries a $\mathsf K$-action.  We will determine the $\mathsf K$-virtual representation $\chi(K_v, \Theta_w)$ explicitly, under assumption $({\bf A.1})$ for the vector $v$. 

\vskip.1in

\addtocounter{theorem}{0}
\begin {theorem}\label{t1} Consider $\zeta=(x,y)\in A[a]\times \widehat A[b]$ of order $\delta$. Then the trace of $\zeta$ is given by $$\text {Trace }(\zeta, \chi(K_v, \Theta_w))=\frac{d_v^2}{d_v+d_w}\binom{d_v/\delta+d_w/\delta}{d_v/\delta}.$$ 
\end {theorem} 
\vskip.1in

\proof 
For future use, we prove a slightly more general version of the Theorem, valid in the context of Bridgeland stability conditions. We work with moduli spaces $K_v(\sigma)$ of Bridgeland stable objects, for certain generic Bridgeland stability conditions $\sigma$. As it is well-known, the Gieseker moduli spaces emerge as particular cases. In fact, the spaces $K_v(\sigma)$ are all smooth birational models of the Gieseker space $K_v$, see \cite {Y3}, Corollary 3.33. The objects parametrized by $K_v(\sigma)$ are certain $2$-step complexes $E^{\bullet}$ in the derived category of $A$. These complexes have cohomology sheaves only in degree $-1$ and $0$ $$H^{-1} (E^{\bullet})\in \mathcal F \,\,\,\text{ and }\,\,\,\, H^0(E^{\bullet})\in \mathcal T,$$ for a certain torsion pair $(\mathcal F, \mathcal T)$; in particular $\text{Hom}(\mathcal T, \mathcal F)=0$. While it is possible to give a more precise description of these complexes, this will not be needed for our argument. 

The discussion of the previous subsections carry over in this context as well. In particular, the moduli space $K_v(\sigma)$ carries an action of $\zeta\in A[a]\times \widehat A[b]$: $$\zeta\cdot E^{\bullet} =t_{-x} ^{\star} E^{\bullet} \otimes y^{-1}.$$ The action similarly lifts to the theta bundle and to the signed sum of cohomologies. We will show that $$\text{Trace }(\zeta, \chi(K_v(\sigma), \Theta_w))=\frac{d_v^2}{d_v+d_w}\binom{d_v/\delta+d_w/\delta}{d_v/\delta}.$$ The proof is an application of the Lefschetz-Riemann-Roch theorem, as in \cite {wirtinger}. The details, both quantitative and qualitative, are however different. \vskip.1in

{\it Step 1.} We first find the fixed points of the action of $\zeta$:$$t_{-x}^{\star} E^{\bullet}\cong E^{\bullet} \otimes y^{-1}, \,\,\, E^{\bullet} \in K_v(\sigma).$$

Possibly replacing $a$ and $b$ by some of their divisors (thus without insisting on $a$ and $b$ being the gcd of the relevant numerics), we may assume that the order of $x$ is $a$ and the order of $y$ is $b$. Of course the following still hold $$a\,|\,\chi, \,\,b\,|\,r \implies \gcd (a,b)=1, \text{ and } ab|k.$$ We have $\delta=ab$. 

Now, $$t_{-x}^{\star} E^{\bullet}\cong E^{\bullet}\otimes y^{-1}\implies t_{ax}^{\star} E^{\bullet}\cong E^{\bullet}\otimes y^{a}\implies E^{\bullet}\cong E^{\bullet}\otimes y^a.$$ Since the order of $y$ is coprime to $a$, we obtain $$E^{\bullet}\cong E^{\bullet}\otimes y \text { and consequently } t_{x}^{\star} E^{\bullet}\cong E^{\bullet}.$$
 
Consider the abelian cover $$p:A\to A''=A/\langle x\rangle.$$ The Galois group $G$ of the cover is generated by translations by $x$. Let $\widehat G$ be the dual group, and pick a generator $x''$ of $\widehat G$. This corresponds to a line bundle $x''\to A''$ of order $a$ which determines the cover $p$. Then $$E^{\bullet}=p^{\star} E^{\bullet ''}$$ for some complex $E^{\bullet''}$ on $A''$ with $$\text{ rank } E^{\bullet''}=r,\,\,\, \chi(E^{\bullet''})=\frac{\chi}{a}.$$ We have $$p^{\star} E^{\bullet''}=E^{\bullet}\cong E^{\bullet}\otimes y=p^{\star} (E^{\bullet''}\otimes y'')$$ where $$y=p^{\star} y''.$$ Pushing forward under $p$, we obtain $$E^{\bullet''}=E^{\bullet''}\otimes y''\otimes x''^{\ell}$$ for some $\ell$. Replacing $y''$ by $y''\otimes x''^{\ell}$ we may assume $$E^{\bullet''}=E^{\bullet''}\otimes y''.$$ Note that $$y^{b}=p^{\star} y''^{b}=1\implies y''^{b}\in \langle x''\rangle \implies y''^{ab}=1.$$ Therefore, the order $\Delta$ of $y''$ satisfies $\Delta|ab$. We will show shortly that in fact $\Delta=b$. 

Now, let $$\pi:A'\to A''$$ be the cover determined by $y''$ which has degree $\Delta$. Let $G'$ denote the Galois group of the cover, and let $\widehat G'$ be the dual group which is generated by $y''$. We collect the following facts about complexes over arbitrary abelian covers:

\begin {itemize}
\item [(i)] if $E^{\bullet''}\cong E^{\bullet''}\otimes y'',$ then $E^{\bullet''}$ is the pushforward of a complex $E^{\bullet '}\to A'$:$$E^{\bullet ''}=\pi_{\star} E^{\bullet'}$$ 
\item [(ii)] $\pi_{\star} E^{\bullet'}_1=\pi_{\star} E^{\bullet'}_2$ iff $E^{\bullet '}_1=\beta^{\star} E^{\bullet '}_2$ for some $\beta\in G'$ 
\item [(iii)] $$\pi^{\star} \pi_{\star} E^{\bullet '}=\bigoplus_{\beta\in G'} \beta^{\star} E^{\bullet'}$$ 
\item [(iv)]  the action of $y''\in \widehat G'$ on $\pi^{\star}\pi_{\star} E^{\bullet '}$ leaves each of the factors $\beta^{\star}E^{\bullet '}$ invariant. The weight of the action on each factor equals $\langle y'', \beta\rangle.$
\end {itemize} 
Specifically, (i) is proven in Proposition 2.5 of \cite {BM}, while (ii), (iii) and (iv) follow by imitating the proofs of Lemmas $2.1$ and $2.2$ in \cite {NR}. 

We conclude $$E^{\bullet ''}=\pi_{\star} E^{\bullet '}$$ where $$\text{ rank } E^{\bullet '}=\frac{r}{\Delta}, \, \chi(E^{\bullet '})=\chi(E^{\bullet ''})=\frac{\chi}{a}.$$ In particular, $\Delta|r$ hence $\Delta$ is coprime to $a$. Since $\Delta$ divides $ab$, we obtain that $\Delta$ divides $b$. On the other hand, $$y''^{\Delta}=1\implies y^{\Delta}=p^{\star} y''^{\Delta}=1\implies b|\Delta.$$ Hence $\Delta =b$. Thus all $\zeta$-fixed complexes are of the form $$E^{\bullet}=p^{\star} \pi_{\star} E^{\bullet '}, \text{ rank }E^{\bullet '}=\frac{r}{b}, \, \chi(E^{\bullet '})=\frac{\chi}{a}.$$ 

It is clear that the map $p$ is uniquely determined by $\zeta=(x, y)$. The same is true about $\pi$. To see this, note that the preimage of $y$ under $p$ contains at most one element of order $b$ since the preimage of $0$ contains at most one such element. Indeed, $p^{\star}(0)$ consists only in powers of $x''$ which all have order dividing $a$. This fixes $y''$ and therefore $\pi$. 

At this point we briefly discuss stability. Consider the stability condition over $A'$ $$\sigma'=\pi^{\star}\sigma'',$$ where $\sigma''$ over $A''$ is chosen so that $$p^{\star}\sigma''=\sigma.$$ It can be immediately checked from definitions that $E^{\bullet}$ is $\sigma$-stable implies that $E^{\bullet'}$ is $\sigma'$-stable, and conversely the stability of $E^{\bullet'}$ implies (semi)stability of $E^{\bullet}$. These choices of stability conditions will be assumed below, and will be suppressed from the notation. 

To describe the $\zeta$-fixed loci over $K_v(A)$, we fix the determinant and determinant of Fourier Mukai of $E^{\bullet}$. We claim that there are $(ab)^3$ fixed loci all isomorphic to $K_{v'}(A')$ for $$v'=\left ( \frac{r}{b}, 0, \frac{\chi}{a}\right).$$ We first calculate the determinant $$\det E^{\bullet}=p^{\star} \det \pi_{\star} E^{\bullet '}=\mathcal O\implies \det \pi_{\star} E^{\bullet'}\in \langle x'' \rangle.$$ Writing $$\det E^{\bullet'}=\pi^{\star} M''$$ we obtain by the proof of Lemma 3.4 in \cite {NR} that $$\det \pi_{\star} E^{\bullet'}=\det \pi_{\star} \pi^{\star} M''=\det \left(\bigoplus_{j=1}^{b} M''\otimes y''^{j}\right)=M''^{b}\in \langle x'' \rangle \implies M''\in \langle x'' \rangle+\widehat A''[b],$$ where above we used that the order of $x''$ is $a$ which is coprime to $b$. Thus $$\det E^{\bullet'}\in\pi^{\star} \langle x''\rangle+\pi^{\star} \widehat A''[b].$$ There are $ab^3$ choices for the determinant of $E^{\bullet'}$, since the torsion point $y''\in \widehat A''[b]$ pulls back trivially to $A'$. The calculation of the Fourier-Mukai determinant is similar. (In the calculation below, for simplicity of notation, we denote by $\,\widehat{}\,$ the Fourier-Mukai over the appropriate abelian varieties, without assuming that these transforms are represented by sheaves.) We have 
$$\det \widehat {E^{\bullet}}=\det \widehat p_{\star}\, \widehat \pi^{\star}\, \widehat {E^{\bullet'}}=\mathcal O.$$ Write $$\det \left(\widehat \pi^{\star} \widehat {E^{\bullet'}}\right)=\widehat p^{\star} M$$ and observe that this gives $$\det \widehat p_{\star} \,\widehat p^{\star} M=M^{a}=\mathcal O\implies M\in A[a].$$ This gives $a^3$ options for $\widehat p^{\star} M$ since $x\in A[a]$ pulls back trivially. Finally, this identifies $\det \widehat E^{\bullet'}$ up to the $b$ elements in the kernel of $\pi$: $$\pi(\det {\widehat E^{\bullet'}})=p(M)\in p(A[a])\implies \det {\widehat {E^{\bullet'}}}\in \pi^{-1} (p(A[a])).$$

We obtain $(ab^3)(a^3b)$ fixed loci, one for each choice of determinant and determinant of Fourier-Mukai of $E^{\bullet'}$. However, this answer does not account for repetitions. We observe that $$p^{\star} \pi_{\star} E^{\bullet'}_1=p^{\star} \pi_{\star} E_2^{\bullet'}\iff \pi_{\star} E_1^{\bullet'}=\pi_{\star} E^{\bullet'}_2\otimes x''^{\alpha} \iff E^{\bullet'}_1=\beta^{\star} E^{\bullet'}_2\otimes \pi^{\star}x''^{\alpha},\,\,\, \beta\in G'.$$ The second statement above follows pushing forward via $p$ and using stability. The third statement is contained in (ii). As a result, we are left with only $(ab)^3$ fixed loci. Indeed, via the above equivalence we only obtain $b^3$ choices for the determinant, if we require $$\det E^{\bullet'}\in \pi^{\star} \widehat A''[b].$$ Similarly, the determinant of the Fourier-Mukai dual can be fixed in $a^3$ ways. It is easy to see there are no other repetitions. This yields $(ab)^3$ distinct fixed loci.
\vskip.1in

{\it Step 2.} We apply the Lefschetz-Riemann-Roch theorem to calculate the trace in the theorem by summing the contributions from the fixed loci $K_{v'}(A')$. We obtain $$\text {Trace} (\zeta, \chi(K_v, \Theta_w))=(ab)^3\int_{K_{v'}(A')}\text {Todd }(K_{v'}(A'))\,\cdot i^{\star} \text{ch}_{Z}(\Theta_w) (\zeta)\,\cdot \prod_{z\neq 1} \left(\text{ch}_{-\langle \zeta^{-1}, z\rangle} N_{z}^{\vee}\right)^{-1}.$$ Here $Z$ is the cyclic group generated by $\zeta$. The notation we used is as follows: $$i^{\star} \text{ch}_{Z} (\Theta_w)\in H_{Z}^{\star} (K_{v'}(A'))=\text{Rep}(Z)\otimes H^{\star}(K_{v'}(A'))$$ is the restriction of the equivariant Chern character via inclusion $$i:K_{v'}(A')\to K_v(A),$$ followed by the evaluation against $\zeta\in Z$ in the representation ring. The normal bundle of the inclusion $i$ splits into eigenbundles $N_z$ indexed by elements in the character group $z\in \widehat Z$. Finally, we write $$\text{ch}_t(N)=\prod_{i}(1+t e^{x_i})$$ for any bundle $N$ with Chern roots $x_1, \ldots, x_{\ell}$. The prefactor $(ab)^3$ comes from the fact that all fixed loci will have identical trace contributions. \vskip.1in

{\it Step 3.} We evaluate the integral above explicitly. We begin by computing the normal bundles $N_z$ in the expression above. We will compute the eigenvalues of the action of $\zeta$ on $T_{E^{\bullet}} K_v$. We first consider the similar action on $T_{E^{\bullet}} \uv$ and use the morphism $$\alpha:\uv\to A\times \widehat A$$ to find the eigenvalues on the fiber.

The tangent space to $\uv$ at a fixed point $E^{\bullet}$ was calculated in \cite {I}: 
\begin{eqnarray*}T_{E^{\bullet}} \uv&=&\text {Ext}^1(E^{\bullet}, E^{\bullet})=\text{Ext}^1 (p^{\star} \pi_{\star} E^{\bullet'}, p^{\star} \pi_{\star} E^{\bullet'})=\text{Ext}^1(\pi_{\star} E^{\bullet'}, \pi_{\star} E^{\bullet'}\otimes p_{\star} \mathcal O)\\&=&\bigoplus_{\alpha=0}^{a-1}\text{Ext}^1 (\pi_{\star} E^{\bullet'}, \pi_{\star} E^{\bullet'}\otimes x''^{\alpha}) =\bigoplus_{\alpha=0}^{a-1} \text{Ext}^1(\pi_{\star} E^{\bullet'}, \pi_{\star} (E^{\bullet'}\otimes \pi^{\star} x''^{\alpha}))\\&=&\bigoplus_{\alpha=0}^{a-1}\text{Ext}^1(\pi^{\star}\pi_{\star}E^{\bullet'}, E^{\bullet'}\otimes \pi^{\star} x''^{\alpha})=\bigoplus_{\alpha=0}^{a-1} \bigoplus_{\beta\in G'} \text{Ext}^1(\beta^{\star} E^{\bullet'}, E^{\bullet'}\otimes \pi^{\star} x''^{\alpha}).\end{eqnarray*}We claim that $$T_{\alpha, \beta}=\text{Ext}^1(\beta^{\star} E^{\bullet'}, E^{\bullet'}\otimes x''^{\alpha})$$ are the isotypical components of the tangent space. Indeed, by (iv) in {\it Step 1}, $y$ acts on each summand via the root of unity $\langle \beta, y''\rangle$ of order $b$. The action of $x$ also leaves the subspace invariant, and the action has weight $\langle x, x''^{\alpha}\rangle$ which is root of unity of order $a$. As $\alpha, \beta$ vary, we obtain all the $(ab)$-roots of unity as eigenvalues. 

We calculate the Chern roots of the eigenbundles $T_{\alpha, \beta}$. Clearly, by Hirzebruch-Riemann-Roch, the Chern character of the virtual bundles $$\sum_{i=0}^{2} (-1)^i \text {Ext}^i(\beta^{\star} E^{\bullet'}, E^{\bullet'}\otimes \pi^{\star} x''^{\alpha})$$ must stay constant as $\beta$ and $\pi^{\star} x''^{\alpha}$ vary in the abelian varieties $A'$ and $\widehat A'$. (The index $i$ is checked to have the correct range. Indeed, vanishing of the ext's for indices $i\leq -1$ is built in the definition of Bridgeland stabilities, while vanishing for $i\geq 3$ follows by duality.) 
Furthermore, for $\alpha, \beta$ not both trivial we have $$\text{Ext}^0(\beta^{\star} E^{\bullet'}, E^{\bullet'}\otimes \pi^{\star} x''^{\alpha})=\text{Ext}^2(\beta^{\star} E^{\bullet'}, E^{\bullet'}\otimes \pi^{\star} x''^{\alpha})=0$$ while for trivial $\alpha$ and $\beta$ the two dimensions are $1$, by stability of $E^{\bullet'}$. Indeed, for nontrivial $(\alpha, \beta)$, we calculate by duality $$\text{Ext}^2(\beta^{\star} E^{\bullet'}, E^{\bullet'}\otimes \pi^{\star} x''^{\alpha})=\text{Ext}^0(E^{\bullet'}\otimes \pi^{\star} x''^{\alpha}, \beta^{\star} E^{\bullet'})=\text{Ext}^0((\beta^{-1})^{\star} E^{\bullet'}, E^{\bullet'}\otimes \pi^{\star} x''^{-\alpha})$$ so it suffices to prove the statement about $\text{Ext}^0$. This is immediate since any non-zero morphism $$\beta^{\star} E^{\bullet'}\to E^{\bullet'}\otimes \pi^{\star} x''^{\alpha}$$ is an isomorphism inducing a map $$\pi_{\star} \beta^{\star} E^{\bullet'}=\pi_{\star}E^{\bullet'}\to \pi_{\star} E^{\bullet'}\otimes x''^{\alpha}.$$ Comparing determinants, we must have $x''^{\alpha r}=0 \implies x''^{\alpha}=0$ since $(r, a)=1$. Therefore, $\alpha=0$. Now, using the isomorphism $$\beta^{\star} E^{\bullet'}\to E^{\bullet'}$$ and letting $$q:A'\to A'/\langle \beta \rangle$$ be the projection determined by $\beta$, we obtain that $E^{\bullet'}$ is a pullback from the quotient. Evaluating Euler characteristics, we obtain that $\chi(E^{\bullet'})=\frac{\chi}{a}$ is divisible by the order of $\beta$. But $\text{ord} (\beta)|b$ and $(b, \chi)=1$, hence $\text{ord} (\beta)=1$ and $\beta=1$. 

We are now in the position to calculate the normal bundles $N_z$. The isotypical components correspond to nontrivial pairs $(\alpha, \beta)$. Each isotypical component has dimension $2\frac{r\chi}{ab}:=\ell+2.$ We just argued the Chern roots of $T_{\alpha, \beta}$ equal $$0, 0, x_1, \ldots, x_{\ell},$$ with $x_i$ the roots of the tangent bundle of $K_{v'}(A')$. This follows by comparison with $\alpha, \beta$ trivial; the two trivial roots come from the four dimensional base of the Albanese map, after canceling two trivial factors corresponding to infinitesimal automorphisms and obstructions over $\mathfrak M_{v'}(A')$.
 
\vskip.1in
{\it Step 4.} 
With this understood, we calculate $$\left(\prod_{z\neq 1}\text{ch}_{-\langle \zeta^{-1}, z\rangle} N_z^{\vee}\right)^{-1}=\prod_{\xi\neq 1}\left((1-\xi)^{2} \prod_{i=1}^{\ell}\left(1-\xi e^{-x_i}\right)\right)^{-1}=\frac{1}{\delta^2}\prod_{i=1}^{\ell} \frac{1-e^{-x_i}}{1-e^{-\delta x_i}},$$ where $\xi=\langle \zeta^{-1}, z\rangle$ runs through the non-trivial $\delta$-roots of $1$. 

We now claim that \begin{equation}\label{izp}i^{\star}_Z \Theta_w\cong \Theta^\delta_{w'},\end{equation} where the last bundle carries a trivial $Z$-linearization and the vector $w'$ is specified below. Indeed, let $\Theta''\to A''$ be the symmetric polarization such that $$p^{\star} \Theta''=\Theta^a \implies \chi (\Theta'')=a.$$ Also, write $\Theta'=\pi^{\star} \Theta'',$ hence $\chi(\Theta')=ab$. Let $w=(r', k\Theta, \chi')$ where $r'=rh, \,\,\chi'=-\chi h$. We introduce the Mukai vectors $$w''=\left(r', \frac{k}{a} \Theta'', \frac{\chi'}{a}\right) \text { over } A'', \text { and } $$ $$w'=\left(\frac{r'}{b}, \frac{k}{ab} \Theta', \frac{\chi'}{a}\right)\text { over } A'.$$ Observe that $$p_{\star} w=aw'', \,\, \pi^{\star} w''=bw'.$$ These equalities imply that \begin{equation}\label{iotap} i_{1}^{\star} \Theta_w=\Theta_{w''}^a, \,\, i_2^{\star} \Theta_{w''}=\Theta_{w'}^b,\end{equation} which together then give \eqref{izp}. Here, we factored $$i:\mathfrak M_{v'}(A')\to \mathfrak M_v(A), \,\,\, i=i_1\circ i_2$$ where $$i_1: \mathfrak M_{v''}(A'')\to \mathfrak M_v(A),\,\,\,\,\, i_2:\mathfrak M_{v'}(A')\to \mathfrak M_{v''}(A'')$$ are induced by pullback by $p$ and pushforward by $\pi$ respectively, and the corresponding Mukai vectors are $$v'=\left(\frac{r}{b}, 0, \frac{\chi}{a}\right),\,\,\ v''=\left(r, 0, \frac{\chi}{a}\right).$$ For instance, to justify the first identity in \eqref{iotap}, note that tautologically we have $$\iota_1^{\star} \Theta_{F}=\Theta_{p_{\star}F},$$ for any complex $F$ representing the vector $w$. Equation \eqref{iotap} follows by recalling the normalization conventions for theta bundles in Section \ref{sec2}. Specifically, if $F$ is suitably normalized as in \eqref{conv1}, then $p_{\star}F$ represents the vector $aw''$ and also satisfies convention \eqref{conv1}.  \vskip.05in

The statement about the $Z$-action in \eqref{izp} will be proved in {\it Step 5} below. 

\vskip.1in

Substituting into Lefschetz-Riemann-Roch, we find 
\begin{eqnarray*} \text{Trace} (\zeta, \chi(K_v, \Theta_w))&=&\delta^3\int_{K_{v'}(A')} \prod_{i=1}^{\ell} \frac{x_i}{1-e^{-x_i}} \cdot \frac{1}{\delta^2}\prod_{i=1}^{\ell} \frac{1-e^{-x_i}}{1-e^{-\delta x_i}} \cdot \text {ch} (\Theta_{w'}^\delta)\\ &=& \delta \int_{K_{v'}(A')}  \prod_{i=1}^{\ell} \frac{x_i}{1-e^{-\delta x_i}} \cdot \text {ch}(\Theta_{w'}^\delta)\\ & =& \delta \chi(K_{v'}, \Theta_{w'}) = \delta \frac{d_{v'}^2}{d_{v'}+d_{w'}} \binom {d_{v'}+d_{w'}}{d_{v'}}. \end {eqnarray*} The last line follows from the backward application of Hirzebruch-Riemann-Roch and \eqref{eulerc}. The proof is completed by observing that $$d_{v'}=\frac{d_v}{\delta}\text { and } d_{w'}=\frac{d_w}{\delta}.$$\vskip.1in

{\it Step 5.} We explain now that the action of $\zeta$ in the fiber of $\Theta_w\to K_v$ over each $\zeta$-fixed point is trivial, as claimed in {\it Step 4}. The idea of the proof is similar to that of Remark $1$ in \cite {wirtinger}. Since the details are different, we include the argument for completeness. Let $E^{\bullet}=p^{\star}\pi_{\star}E^{\bullet'}$ be a $\zeta$-fixed point of $K_v(A)$. Consider $$t:A\times \widehat A\to \mathfrak M_v(A),\,\,\, t(\lambda, \mu)=t_\lambda^{\star}E^{\bullet}\otimes \mu$$ the restriction of the morphism $\tau$ to $\{E^{\bullet}\}\times A\times \widehat A.$ By the construction in Subsection \ref{ga}, it suffices to explain that the identification \begin{equation}\label{zetaequiv}t^{\star}\Theta_w\simeq (a, b)^{\star} \left(\Theta^{\frac{-\chi k}{a^2}}\boxtimes \widehat \Theta^{-\frac{rk}{b^2}}\otimes \mathcal P^{\,h\cdot \frac{r\chi}{ab}}\right)\end{equation} obtained by restricting \eqref{pback2} to $\{E^{\bullet}\}\times A\times \widehat A$ is $\zeta$-equivariant. 

To this end, consider the fiber diagram  
$$\xymatrix{A^+\ar[r]^{\,\,\,\,\pi^+} \ar[d]^{p^+} & A \ar[d]^{p} \\ A'\ar[r]^{\,\,\,\,\pi} & A''}.$$ The constructions in {\it Step 1} show that $\pi^+:A^+\to A$ is the \'etale cover determined by $p^{\star}y''=y$, so that $$\widehat {\pi^+}:\widehat A\to \widehat {A^+}$$ is the quotient by the translation action by $y$. The natural map $$(p, \widehat {\pi^+}):A\times \widehat A\to A''\times \widehat {A^+}$$ is the quotient by the action of $\zeta=(x, y)$. Since $t$ is $\zeta$-invariant, we obtain a morphism $$\overline t:A''\times \widehat {A^+}\to \mathfrak M_v(A),\,\,\,\,\, \overline t \circ (p, \widehat {\pi^+})=t.$$ Let $$N:A''\to A,\,\,\, \widehat N^+:\widehat {A^+}\to \widehat A$$ denote the two norm maps corresponding to the morphisms $p$ and $\widehat {\pi^+}$ respectively. Then $$N\circ p=a,\,\,\, \widehat {N^+}\circ \widehat {\pi^+}=b\implies (N, \widehat {N^+})\circ(p, \widehat {\pi^+})=(a, b).$$ To establish that \eqref{zetaequiv} holds $\zeta$-equivariantly, we first factor out the action of $\zeta$, and prove that over the quotient $A''\times \widehat {A^+}$ we have \begin{equation}\label{tprime}{\overline t}^{\,\star}\Theta_w\simeq (N, \widehat {N^+})^{\star} \left(\Theta^{\frac{-\chi k}{a^2}}\boxtimes \widehat \Theta^{-\frac{rk}{b^2}}\otimes \mathcal P^{\,h\cdot \frac{r\chi}{ab}}\right).\end{equation} 

The isomorphism \eqref{tprime} will be shown using the see-saw principle. We verify it over $A''\times \{\mu^+\}$ for all $\mu^+\in \widehat {A^+}$. The restriction to $\{\lambda''\}\times \widehat {A^+}$ is similar and will be omitted. We show \begin{equation}\label{tdouble}(t'')^{\star}\Theta_w\simeq N^{\star}\left(\Theta^{-\frac{\chi k}{a^2}}\otimes \widehat {N^+}(\mu^+)^{\,h\cdot \frac{r\chi}{ab}}\right),\end{equation} where $$t'':A''\to \mathfrak M_v(A) \text{ is the restriction of } \overline t,\text { that is } t''(\lambda'')={\overline t}(\lambda'', \mu^+).$$ We first determine the morphism $t''$. Write $$\mu^+=q^{\star}\mu''$$ for some $\mu''\in \widehat A''$, where $q=p\circ \pi^+:A^+\to A''.$ It follows from the definitions that $$\overline t(\lambda'', \mu^+)=t(\lambda, \mu)=t_{\lambda}^{\star}E^{\bullet}\otimes \mu$$ whenever $$p(\lambda)=\lambda'',\,\,\,\, \mu^+=(\pi^+)^{\star}\mu.$$ In our case, we can take $\mu=p^{\star}\mu''$, so that $$t''(\lambda'')=\overline t(\lambda'', \mu^+)=t_{\lambda}^{\star}\,E^{\bullet}\otimes p^{\star}\mu''=t_{\lambda}^{\star}\, p^{\star}\,E^{\bullet''}\otimes p^{\star} \mu''=p^{\star}\left(t_{\lambda''}^{\star} E''\otimes \mu''\right).$$ Here, we used that the $\zeta$-fixed points take the form $$E^{\bullet}=p^{\star} E^{\bullet''} \text {with } E^{\bullet ''}=\pi_{\star}E'^{\bullet}.$$ In conclusion, 
$$t''=\iota_1\circ \tau''$$ where $$\iota_1:\mathfrak M_{v''}(A'')\to \mathfrak M_v(A)$$ is the morphism induced by the pullback $p^{\star}$ already encountered in {\it Step 4}, and $$\tau'':A''\to \mathfrak M_{v''}(A''),\,\, \tau''(\lambda'')=t_{\lambda''}^{\star}E^{\bullet ''}\otimes \mu''$$ is the translation map. 

With this understood, we compute the left hand side of \eqref{tdouble} making use of equation \eqref{iotap}: $$(t'')^{\star}\, \Theta_w=(\iota_1\circ \tau'')^{\star}\,\Theta_w=(\tau'')^{\star}\, (\iota_1)^{\star}\, \Theta_w=(\tau'')^{\star}\,\Theta_{w''}^a.$$ The pullback of the theta bundle $\Theta_{w''}$ under $\tau''$, all the way to the product $A''\times \widehat {A''}$, was determined in Lemma \ref{lemmal}. (We remarked above that the Lemma also holds for nonprincipal polarizations.) Restricting to $A''\times \{\mu''\}$, we find $$(\tau'')^{\star} \Theta_{w''}=\Theta''^{-\frac{\chi k}{a^2}}\otimes \left(\mu''\right)^{h \frac{r\chi}{a}}\implies (t'')^{\star} \Theta_w=\Theta''^{-\frac{\chi k}{a}}\otimes \left(\mu''\right)^{h r\chi}.$$
For the right hand side of \eqref{tdouble}, we use the standard identities $$N^{\star}\Theta=\Theta''^a, \,\,\, N^{\star} \widehat {N^+}q^{\star} (\mu'')=(\mu'')^{ab}.$$ Hence, $$N^{\star}\left(\Theta^{-\frac{\chi k}{a^2}}\otimes \widehat {N^+}(\mu^+)^{\,h\cdot \frac{r\chi}{ab}}\right)=\Theta''^{-\frac{\chi k}{a}}\otimes \left(\mu''\right)^{h r\chi}.$$ This establishes \eqref{tdouble}, and completes the argument. 
\qed

\subsection{The calculation of the Verlinde bundles} \label{vbb}We are now in the position to determine the Verlinde bundles in degree $0$. Assuming ${\bf (A.1)}$ and ${\bf (A.2)}$, 
we prove

\begin {theorem} \label{thm21}We have $$\mathbf E(v, w)=\bigoplus_{\zeta} \left(\mathbf W_{-\frac{\chi}{a}, \frac{k}{a}}\boxtimes \mathbf W^{\dagger}_{\frac{r}{b}, -\frac {k}{b}}\otimes \mathcal P^{-h}\right)\otimes  \ell_{\zeta}^{\,\oplus\mathsf m_{\zeta}}.$$ The sum is indexed by torsion line bundles $\zeta\to A\times \widehat A$ of order dividing $(a, b).$ An element $\zeta$ of order exactly $\omega$ comes with multiplicity $$\mathsf m_{\zeta}=\frac{1}{d_v+d_w} \sum_{\delta|ab} \frac {\delta^4}{(ab)^2} \left\{\frac{ab/\omega}{\delta}\right\} \binom{d_v/\delta+d_w/\delta}{d_v/\delta}.$$
\end {theorem} 

Recall that in the above, for each line bundle $\zeta$ of order $(a,b)$ over $A\times \widehat A$, we fix {\it one} root $\ell_{\zeta}$ such that $$\left(-\frac{\chi}{a}, \frac{r}{b}\right)\ell_{\zeta}=\zeta\implies (-\chi, r)\ell_{\zeta}=0.$$ Each $\ell_{\zeta}$ corresponds to a character of $A[-\chi]\times \widehat A[r]$ which is uniquely defined only up to characters of $A[-\chi/a]\times \widehat A[r/b].$

\begin{remark} A more general class of semihomogeneous vector bundles $$\mathbf W(\mathsf P)\to A\times \widehat A,$$ depending on a triple $\mathsf P$ of rational numbers, is constructed and studied in \cite{class}. For instance, for a triple written in lowest terms $$\mathsf P=\left(\frac{b}{a}, \frac{d}{c}, h\right),$$ where $(a, c)$ are odd positive and $h\in \mathbb Z$, we have $$\mathbf W(\mathsf P)=\mathbf W_{a,b}\boxtimes \mathbf W_{c,d}^{\dagger}\otimes \mathcal P^h.$$ For general triples $\mathsf P$, the bundles $\mathbf W(\mathsf P)$ do not admit such simple expressions.  Conjecturally, for arbitrary numerics, the Verlinde bundle can be written in terms of the irreducible building blocks $\mathbf W(\mathsf P)$ and torsion points with explicit multiplicities, in a fashion compatible with the Fourier-Mukai symmetry of Theorem \ref{tttt}. This will be investigated in more detail in \cite{class}. 
\end{remark}

\proof The proof of the theorem follows the strategy laid out in \cite {wirtinger}, with a few modifications. We give the relevant details here. 

Let us assume first that $\Theta_w\to K_v$ carries no higher cohomology. We noted in \eqref {pback1} that $$(-\chi, r)^{\star} \mathbf E(v, w)=H^0(K_v, \Theta_w)\otimes (\Theta^{-\chi}\boxtimes \widehat \Theta^{-r})^{k} \otimes \left((-\chi, r)^{\star} \mathcal P^{-h}\right).$$ This identifies $\mathbf E(v, w)$ up to $(-\chi, r)$-torsion line bundles. In Section $4.0.2$ of \cite {wirtinger} we observed the following equivariant identifications $$(-\chi)^{\star} {\mathbf W}_{-\frac{\chi}{a}, \frac{k}{a}}=\left(\Theta^{-\chi}\right)^ k \boxtimes R_1$$ and $$ r^{\star} {\mathbf W}^{\dagger}_{\frac{r}{b}, -\frac{k}{b}}=\left(\widehat \Theta^{r}\right)^{-k}\boxtimes R_2$$ for $R_1$ a representation of ${\mathsf H}[-\chi]$ of dimension $(\chi/a)^2$ and central weight $-k$, while $R_2$ is a representation of $\widehat {\mathsf H} [r]$ of dimension $(r/b)^2$ and central weight $k$. Therefore, $$(-\chi, r)^{\star} \left({\mathbf W}_{-\frac{\chi}{a}, \frac{k}{a}}\boxtimes {\mathbf W}_{\frac{r}{b}, -\frac{k}{b}}^{\dagger}\right)=\left(\Theta^{-\chi}\boxtimes \widehat \Theta^{-r}\right)^{k} \boxtimes R$$ where $$R=R_1\boxtimes R_2$$ is the product representation of ${\mathsf H}[-\chi]\times {\widehat {\mathsf H}} [r].$ It suffices to explain that $\mathsf H[-\chi]\times \widehat {\mathsf H}[r]$- equivariantly we have \begin {equation}\label{thm3proof} H^0(K_v, \Theta_w)=R\otimes \bigoplus_{\zeta} \ell_{\zeta}^{\,\oplus \mathsf m_{\zeta}}\end{equation} where $\ell_{\zeta}$ runs over the characters of $A[-\chi]\times \widehat A[r]$ modulo those of $A\left [-\chi/a\right]\times \widehat A[r/b]$.

We make use of the morphism of Theta groups $${\mathsf H}[a]\times \widehat {\mathsf H}[b]\to \mathsf H[-\chi]\times \widehat {\mathsf H}[r]$$ which restricts to $$(\alpha, \beta)\to (\alpha^{-\chi/a}, \beta^{r/b})$$ over the center $$\mathbb C^{\star}\times \mathbb C^{\star}\hookrightarrow \mathsf H[a]\times \widehat {\mathsf H}[b].$$ Furthermore, passing to the finite Heisenberg, two $\widetilde {\mathsf H}[-\chi]\times \widetilde {\mathsf H}[r]$-modules with the central weight $(-k, k)$ are isomorphic if and only if they are isomorphic as representations of the abelian group $A [a]\times \widehat A[b]$, see \cite {wirtinger}. Therefore, it suffices to establish the identification \eqref{thm3proof} equivariantly for the action of $A[a]\times \widehat A[b]$ on both sides.

Crucially, it was explained in Section $3.2$ of \cite {elliptic} that $R_1$ and $R_2$ are the trivial representations of $A[a]$ and $\widehat A[b]$. Same as in \cite {wirtinger}, for $\zeta$ of order exactly $\omega$ dividing $(a, b)$, we use \eqref{thm3proof} to compute \begin{eqnarray*} \mathsf m_{\zeta} &=&\frac{1}{\dim R}\cdot \frac{1}{a^4b^4}\sum_{\pi \in A[a]\times {\widehat A}[b]} \langle \zeta, \pi^{-1}\rangle \text {Trace }\left (\pi, H^0(K_v, \Theta_w)\right)\\ &=&\frac{1}{a^2b^2} \sum_{\delta | ab} \frac{1}{d_v+d_w} \binom{d_v/\delta+d_w/\delta}{d_v/\delta} \left(\sum_{\text{ord }(\pi)=\delta} \langle \zeta, \pi^{-1}\rangle \right)\end {eqnarray*} Theorem \ref{t1} was used here to evaluate the trace. Lemma $4$ of \cite {elliptic} gives the sum $$\sum_{\text{ord }(\pi)=\delta} \langle \zeta, \pi^{-1}\rangle = {\delta^4} \left\{\frac{ab/\omega}{\delta}\right\}.$$ This confirms that the multiplicities $\mathsf m_{\zeta}$ of equation \eqref{thm3proof} agree with the expressions stated in the Theorem. \vskip.1in

We can now remove the assumption on vanishing of higher cohomology. By $({\bf A.2})$, $\Theta_w\to K_v$ belongs to the movable cone, hence it is big and nef on a smooth birational model of $K_v$, by Theorem $7$ in \cite {HT}. Furthermore, the smooth birational models of $K_v$ are obtained as moduli spaces of Bridgeland stable objects $K_v(\sigma)$ for stability conditions $\sigma$ of the type we considered in the proof of Theorem \ref{t1}; see \cite {Y3}, Corollary 3.33 for details. The proof of \cite {Y3} moreover yields the diagram 
$$\xymatrix{\mathfrak M_{v}\ar@{.>}[r]^{\,\,\,\,j} \ar[d]^{\alpha} & \mathfrak M_v(\sigma) \ar[d]^{\alpha^{\sigma}} \\ A\times \widehat A \ar@{=}[r] & A\times \widehat A}$$
where $$\alpha^{\sigma}:\mathfrak M_v(\sigma)\to A\times \widehat A$$ is the Albanese morphism normalized as in equation \eqref{normalb}, and $j$ is birational, regular in codimension $1$, and given by the identity on the common open locus. The two theta line bundles $\Theta_w$ agree under $j$. 

Define $$\mathbf E^{\sigma}(v, w)=(\alpha^{\sigma})_{\star}\Theta_w.$$ The first part of the argument applies verbatim to the moduli of $\sigma$-stable objects. In particular, since $\Theta_w$ carries no higher cohomology over $K_v(\sigma)$, we have $$\mathbf E^{\sigma}(v, w)=\bigoplus_{\zeta} \left(\mathbf W_{-\frac{\chi}{a}, \frac{k}{a}}\boxtimes \mathbf W^{\dagger}_{\frac{r}{b}, -\frac {k}{b}}\otimes \mathcal P^{-h}\right)\otimes \ell_{\zeta}^{\,\oplus\mathsf m_{\zeta}}.$$ The above diagram shows however that $$\mathbf E(v, w)\simeq\mathbf E^{\sigma}(v, w),$$ completing the proof. \qed

\begin {example}\label{example8} {\it Rank $1$.} Let $$v=(1, 0, -n),\,\,\,w=(1, k\Theta, n)$$ with $n$ odd and $k\geq 1$.  Then, by Example \ref{example1}, we have $$\uv\cong A^{[n]}\times \widehat A,\,\,\, \alpha=(-a, 1),$$ and $$\Theta_{w}=(\Theta^{k})_{[n]}\boxtimes \widehat \Theta^{-k}\otimes (a, 1)^{\star} \mathcal P.$$ Therefore, \begin{equation}\label{evw}{\mathbf E}(v,w)=a_{\star} \left((\Theta^{k})_{[n]}\right)\boxtimes \widehat \Theta^{-k} \otimes \mathcal P^{-1}.\end{equation} Theorem \ref{t1} is equivalent to the calculation of the pushforward $$a_{\star}\left( (\Theta^{k})_{[n]}\right)= \bigoplus_{\zeta} \mathbf W_{\frac{n}{a}, \frac{k}{a}}\otimes\ell_{\zeta}^{\,\oplus \mathsf m_{\zeta}},$$ where $\zeta$ are line bundles over $A$ of order $\omega$ dividing $a=\gcd (n, k)$, $\ell_{\zeta}$ is a root of $\zeta$ of order $\frac{n}{a}$, and $$\mathsf m_{\omega}=\frac{1}{k^2}\sum_{\delta|a} \frac{\delta^4}{a^2} \left\{\frac{a/\omega}{\delta}\right\} \binom{k^2/\delta}{n/\delta}.$$ 
To apply the theorem, we invoke the result of Scala who proved the vanishing of higher cohomology of the tautological bundle $(\Theta^k)_{[n]}\to K_{n-1}$, under the assumption that $k\geq 1$, cf. Theorem $5.2.1$ \cite {S}. 

\end{example}

\subsection {Fourier-Mukai symmetries} We can now prove the following Fourier-Mukai comparison, which may be seen as evidence for the strange duality conjecture for abelian surfaces; cf. Conjecture \ref{cj1} of Section \ref{sec2}. Under the assumptions $({\bf A.1})-({\bf A.2})$ made throughout the paper, we establish:

\begin {theorem} \label{tttt} When $c_1(v)$ and $c_1(w)$ are divisible by the ranks $r$ and $r'$, there is an isomorphism $$\widehat{\mathbf E(v, w)}\cong {\mathbf E(w, v)}^{\vee}.$$
\end {theorem}
  
\proof The proof is by direct computation of both sides, using the expression for the Verlinde bundles in Theorem \ref{thm21}. 

To begin, write $$v=(r, r\ell\Theta, \chi), \,\, w=(r', r'\ell'\Theta, \chi').$$ The requirements $\chi(v\cdot w)=0$ and $(r, \chi)=(r', \chi')=1$ give $$2rr'\ell \ell'=-r\chi'-r'\chi\implies r|r' \text { and } r'|r \implies r=r'.$$ While in our setting the numerics impose the restriction $r=r'$ on the ranks, the symmetry in the Theorem is expected to hold true for general numerics satisfying $({\bf A.1})-({\bf A.2})$. In \cite {class} we will offer evidence for this more general statement, but in a partially conjectural setup. 
\vskip.1in

{\it Step 1.} We will reduce to the case $\ell=0$ using twists by line bundles. To this end, we check that the symmetry in the statement of the Theorem is invariant under such twists. Specifically, letting $$v_0=v\exp (- \ell\Theta),\,\,\, w_0=w\exp (\ell\Theta),$$ we consider the two natural  isomorphisms given by tensoring $$i:\uv\to \mathfrak M_{v_0}, \,\, E\mapsto E\otimes \Theta^{-\ell}\text{ and } j: \mathfrak M_w\to \mathfrak M_{w_0}, \,\,\,F\to F\otimes \Theta^{\ell}.$$ It is easy to see that $$i^{\star} \Theta_{w_0}=\Theta_w,\,\, j_{\star}\Theta_v=\Theta_{v_0}.$$ Equation \eqref{deteq} is used here to conclude that tensorization by $\Theta$ does not change the nornalization convention \eqref{conv1} used in the definition of the theta bundles. Next, the Albanese morphisms $$\alpha_0:\mathfrak M_{v_0}\to A\times \widehat A, \,\,\,\beta_0:\mathfrak M_{w_0}\to A\times \widehat A,$$ are respectively given by $$E\mapsto (\det \fm ({E}) , \det E), \,\,\,F\mapsto (\det \fm ({F})\otimes \widehat \Theta^{-r(\ell+\ell')} , \det F \otimes \Theta^{-r(\ell+\ell')}).$$
We claim that \begin{equation}\label{commute}\alpha_0\circ i=\rho \circ \alpha,\,\,\, \beta_0\circ j=\rho'\circ \beta\end{equation} where $$\rho(x, y)=(x-\ell\widehat \Phi(y), y),\,\, \rho'(x, y)=(x+\ell\widehat \Phi(y), y).$$ For instance, the first identity in \eqref{commute} is equivalent to $$\det \fm ({E\otimes \Theta^{-\ell}})=\left(\det \fm ({E})\otimes \widehat \Theta^{-r\ell}\right) \otimes \mathcal P_{-\ell\widehat \Phi(\det E\otimes \Theta^{-r\ell})}.$$ When $E$ is in the Albanese fiber over the origin, so that $$\det E=\Theta^{r\ell} \text{ and }\det \fm (E)=\widehat \Theta^{r\ell},$$ this follows from repeated application of equation \eqref{deteq}. For the general case, we make use of the diagram \eqref{sdia} to write $$E=t_x^{\star}E_0\otimes y,$$ with $E_0$ in the Albanese fiber. The identity above continues to hold by the usual properties of the Fourier-Mukai under tensorization and translation. 

From the first identity in \eqref{commute} we find $$\rho^{\star} \mathbf E(v_0, w_0)=\mathbf E(v, w).$$ In a similar fashion, the second identity in \eqref{commute} gives $$\rho'_{\star}\mathbf E(w, v)=\mathbf E(w_0, v_0)\implies \mathbf E(w, v)=\rho_{\star}\mathbf E(w_0, v_0),$$ using that $\rho'=\rho^{-1}$. 

Assume now that we established the Theorem for the pair $(v_0, w_0)$. The same result will then follow for the pair $(v, w)$. Indeed, it is not hard to see that $\rho=\rho^t$, and hence $$\mathbf E(w, v)^{\vee}=\rho_{\star} \left(\mathbf E(w_0, v_0)^{\vee}\right)=\rho_{\star} \left(\widehat{\mathbf E(v_0, w_0)}\right)=\widehat{{\rho^{\star}\mathbf E(v_0, w_0)}}=\widehat {\mathbf E(v, w)}.$$
\vskip.1in

{\it Step 2.} As a consequence of {\it Step 1}, we may assume $c_1(v)=0$. After clearing primes from the notation, we write $$v=(r, 0, \chi),\, w=(r, r\ell \Theta, -\chi).$$ By assumption, we  have $$\chi=-\frac{d_v}{r} \text{ is an odd negative integer and } \chi+r\ell^2=\frac{d_w}{r}\text{ is an odd positive integer}.$$ To establish \begin{equation}\label{toestablish}\widehat {{\mathbf E}(v, w)}\cong{\mathbf E}(w, v)^{\vee},\end{equation} we explicitly calculate the Verlinde bundles. First, by Theorem \ref{thm2} we have $${\mathbf E}(v, w)=\bigoplus_{\zeta} \left(\mathbf W_{-\frac{\chi}{a}, \frac{r\ell}{a}}\boxtimes \widehat \Theta^{-\ell}\otimes \mathcal P^{-1}\right)\otimes \ell_{\zeta}^{\,\oplus \mathsf m_{\zeta}(v, w)},$$ where $$a=(\chi, r\ell)=(\chi, \ell),\,\,\, \left(-\frac{\chi}{a}, 1\right)\ell_{\zeta}=\zeta.$$ In the sum, $\zeta$ runs over the $(a, 1)$-torsion line bundles over $A\times \widehat A$. Henceforth, we will regard $\zeta$ and $\ell_{\zeta}$ as line bundles pulled back from $A$, without explicitly stating this fact. 

To find ${\mathbf E}(w, v)$, we use calculations similar to those in {\it Step 1} to reduce to degree $0$. In the new notation, we consider the isomorphism $$k: \mathfrak M_w\to \mathfrak M_{w_1},\,\, F\to F\otimes \Theta^{-\ell}, \,\,\, w_1=w\cdot \exp\left({-\ell \Theta}\right), \,\,\,\,v_1=v\cdot \exp\left(\ell\Theta\right).$$ Note that $$w_1=(r, 0, -\chi-r\ell^2), \,\,\,v_1=(r, r\ell \Theta, \chi+r\ell^2).$$ Then,  as before, we have  
$${\mathbf E}(w, v)=\rho^{\star}{\mathbf E}(w_1, v_1).$$ Again by Theorem \ref{thm2}, we have $${\mathbf E}(w_1, v_1)=\bigoplus_{\zeta} \left(\mathbf W_{\frac{\chi+r\ell^2}{a}, \frac{r\ell}{a}}\boxtimes \widehat{\Theta}^{-\ell}\otimes \mathcal P^{-1}\right)\otimes \widetilde{\ell_{\zeta}}^{\,\oplus\mathsf m_{\zeta}(w_1, v_1)}.$$ Here, we used that $(\chi+r\ell^2, r\ell)=(\chi, \ell)=a.$ As above, the sum also runs over the $(a, 1)$-torsion line bundles $\zeta\to A\times \widehat A$, and \begin{equation}\label{tors}\left(\frac{\chi+r\ell^2}{a},1\right)\widetilde\ell_{\zeta}=\zeta.\end{equation}
\vskip.1in
{\it Step 3.} To simplify notation, set $$\mathbf W=\left(\mathbf W_{-\frac{\chi}{a}, \frac{r\ell}{a}}\boxtimes \widehat \Theta^{-\ell}\right)\otimes \mathcal P^{-1},\,\,\,{\mathbf W'}=\left(\mathbf W_{\frac{\Delta}{a}, \frac{r\ell}{a}}\boxtimes \widehat \Theta^{-\ell}\right)\otimes \mathcal P^{-1},$$ where we wrote $$\Delta=\chi+r\ell^2.$$ Then, by the discussion above we have $$\mathbf E(v, w)=\bigoplus_{\zeta} \mathbf W \otimes \ell_{\zeta}^{\,\oplus\mathsf m_{\zeta}(v, w)}, \,\,\, \mathbf E(w, v)=\bigoplus_{\zeta} \rho^{\star} \left(\mathbf W'\otimes \widetilde{\ell_{\zeta}}^{\,\oplus\mathsf m_{\zeta}(w_1, v_1)}\right).$$ It is clear from the specific expressions for the multiplicities ${\mathsf m}_{\zeta}$ given in Theorem \ref{thm2} that these quantities are symmetric in $v$ and $w$, and are invariant under twists, so that $${\mathsf m}_{\zeta}(w_1, v_1)={\mathsf m}_{\zeta}(w, v)={\mathsf m}_{\zeta}(v, w).$$ To complete the proof of \eqref{toestablish}, it suffices to give a correspondence $$\ell_{\zeta}\to \tilde \ell_{\zeta}$$ such that \begin{equation}\label{eqfmzeta}\rho^{\star}\left(\mathbf W'\otimes \widetilde \ell_{\zeta}\right)^{\vee}=\fm\left(\mathbf W\otimes \ell_{\zeta}\right),\end{equation} up to shifts by the index. We first consider the case $\zeta$ trivial, proving that (up to shifts) \begin{equation}\label{is2}\rho^{\star}{\mathbf W'}^{\vee}=\fm ({\mathbf W}).\end{equation} To this end, note first that ${\mathbf W}$ satisfies the index theorem \cite {mukai} with index $0$ if $\ell>0$ or index $4$ if $\ell<0$. This can be checked after pullback. In our case, using \eqref{split}, we find $$\left(-\frac{\chi}{a}, 1\right)^{\star}{\mathbf W}=\left(\Theta^{-\frac{\chi}{a}\cdot \frac{r\ell}{a}}\boxtimes \widehat \Theta^{-\ell}\boxtimes \mathcal P^{\frac{\chi}{a}}\right)\otimes \mathbb C^{\left(\frac{\chi}{a}\right)^2}.$$ Suppose $\ell>0$. The claim about the index follows since the latter line bundle is ample. In turn, this is a consequence of the special form of the Nakai-Moishezon criterion in the context of abelian varieties \cite{bl}. Indeed, a direct calculation as in the last section of \cite {BMOY} shows that the line bundle $$\Theta^{\alpha}\boxtimes \widehat \Theta^{-\beta}\otimes \mathcal P^{\gamma}\to A\times \widehat A$$ is ample whenever $$\alpha>0,\, \beta>0,\, \alpha\cdot \beta-\gamma^2>0.$$ These requirements are satisfied for the numerics we consider. 

We now turn to the proof of \eqref{is2}. We present here a direct argument, referring the reader to the note \cite{class} for similar but more involved statements. We observe first that both bundles in \eqref{is2} are simple and semihomogeneous, as both properties are preserved by pullbacks under isomorphisms and Fourier-Mukai. A direct calculation shows that they  have the same rank $$\text{rank }\mathbf W'=\left(\frac{\Delta}{a}\right)^2,\,\, \text{rank }\fm(\mathbf W)=\chi(\mathbf W)=\left(\frac{\Delta}{a}\right)^2.$$ The slopes of $\rho^{\star}\mathbf W'^{\vee}$ and $\fm(\mathbf W)$ are also directly calculated and seen to match $$\rho^{\star}\left(\Theta^{-\frac{r\ell}{\Delta}}\boxtimes \widehat \Theta^{\ell}\otimes \mathcal P\right)=\Theta^{-\frac{r\ell}{\Delta}}\boxtimes {\widehat\Theta}^{-\frac{\ell\chi}{\Delta}}\otimes \mathcal P^{\frac{\chi}{\Delta}}.$$ In the line above, we used the following identities derived via the see-saw principle $$\rho^{\star} \Theta=\Theta\boxtimes \widehat \Theta^{-\ell^2}\otimes \mathcal P^{\ell},\,\, \rho^{\star}\widehat \Theta=\widehat \Theta,\,\,\, \rho^{\star} \mathcal P=\widehat\Theta^{-2\ell}\otimes \mathcal P.$$ Finally, both bundles in \eqref{is2} are invariant under the involution $(-1, -1)$ over $A\times \widehat A$. The same argument as in Section 2.2 of \cite {wirtinger} shows uniqueness of simple symmetric semihomogeneous bundles with equal rank and determinant, proving \eqref{is2}. 

\vskip.1in

{\it Step 4.} Finally, to establish \eqref{toestablish}, we match the contributions of the torsion points in each irreducible summand \eqref{eqfmzeta}. Fix $\zeta\in \widehat A$ an $a$-torsion line bundle over $A$, viewed as a line bundle over $A\times \widehat A$ by pullback. Let $\ell_{\zeta}\in \widehat A$ be chosen so that $$-\frac{\chi}{a}\ell_{\zeta}=\zeta.$$
The bundle $\mathbf W=\left(\mathbf W_{-\frac{\chi}{a}, \frac{r\ell}{a}}\boxtimes \widehat \Theta^{-\ell}\right)\otimes \mathcal P^{-1}$ is semihomogeneous of rank $\left(-\frac{\chi}{a}\right)^2$ with determinant $$\mathcal D=\left(-\frac{\chi}{a}\right)^2\left(-\frac{r\ell}{\chi}\Theta-\ell\widehat \Theta-\mathcal P\right).$$ By Lemma $6.7$ in \cite {mukai2}, for all $\alpha\in A\times \widehat A$ we have that $$t_{\left(-\frac{\chi}{a}\right)^2\alpha}^{\star}\mathbf W=\mathbf W\otimes \mathcal \phi_{\mathcal D}(\alpha).$$ In the above equation, the line bundle $\mathcal D$ induces over the abelian fourfold $A\times \widehat A$ the Mumford homomorphism $$\phi_{\mathcal D}:A\times \widehat A\to A\times \widehat A,\,\,\,\, (x, y)\mapsto \left(-\frac{\chi}{a}\right)^2\left(-x-\ell \widehat \Phi(y), -\frac{r\ell}{\chi}\Phi(x)-y\right).$$ Pick $y\in \widehat A$ such that $$-\frac{\Delta}{a}\cdot \frac{\chi}{a}y=\ell_{\zeta}$$ and define $$\alpha=(-\ell \widehat \Phi(y), y)\in A\times \widehat A \text{ and } \widetilde{\ell_{\zeta}}=\left(0, \left(\frac{\chi}{a}\right)^2 y\right)\in A\times \widehat A.$$ Our choices guarantee two crucial identities $$\phi_{\mathcal D}(\alpha)=(0, \ell_{\zeta})$$ and $$\rho^{\star}{\widetilde {\ell_{\zeta}}^{\vee}}=-\left(\frac{\chi}{a}\right)^2\alpha.$$ From here, to demonstrate \eqref{eqfmzeta}, we match the contributions of each $\zeta$ to the two sides $$\fm ({\mathbf W}\otimes \ell_{\zeta})=\fm({\mathbf W}\otimes \phi_{\mathcal D}(\alpha))=\fm (t_{\left(\frac{\chi}{a}\right)^2\alpha}^{\star}\mathbf W)=\fm({\mathbf W})\otimes -\left(\frac{\chi}{a}\right)^2\alpha$$ $$=\fm({\mathbf W})\otimes \rho^{\star}{\widetilde {\ell_{\zeta}}^{\vee}}=\rho^{\star}\left({\mathbf W'}\otimes {{\widetilde \ell_{\zeta}}}\right)^{\vee}.$$ To complete the argument, it suffices to note that $$-\frac{\chi}{a}\ell_{\zeta}=\zeta\implies \frac{\Delta}{a}\cdot \left(\frac{\chi}{a}\right)^2y=\zeta\implies \frac{\Delta}{a}\cdot \widetilde{\ell_{\zeta}}=\zeta,$$ as required by \eqref{tors}.
\qed 

\end {document}